\documentclass[a4paper,12pt]{article}

\usepackage[vscale=0.85,includefoot]{geometry}
\usepackage{color}
\usepackage[latin1]{inputenc}
\usepackage[leqno,centertags]{amsmath}
\usepackage{amssymb,latexsym}
\usepackage[english,francais]{babel}
\usepackage{theorem}
\normalbaselines
\def\ifm#1#2{\relax\ifmmode#1\else#2\fi}

\newcommand{\A}{\mathbb{A}}
\newcommand{\R}{\mathbb{R}}
\newcommand{\C}{\mathbb{C}}
\newcommand{\Z}{\mathbb{Z}}
\newcommand{\Q}{\mathbb{Q}}

\newcommand{\M}{\mathbb{M}}

\renewcommand{\P}{\mathbb{P}}

\newcommand{\B}{\mathcal{B}}

\newcommand{\xon}    {\ifm {X_1,\ldots,X_n} {$X_1,\ldots,X_n$}}

\newcommand{\yon}    {\ifm {Y_1,\ldots,Y_n} {$Y_1,\ldots,Y_n$}}

\newcommand{\Fop}    {\ifm {F_1,\ldots,F_p} {$F_1,\ldots,F_p$}}
\newcommand{\iFop}    {\ifm {(F_1,\ldots,F_p)} {$(F_1,\ldots,F_p$)}}

\newcommand{\Qxon}   {\ifm {\Q[X_1,\ldots,X_n]} {$\Q[X_1,\ldots,X_n]$}}

\newcommand{\spar} {\vskip 0.2cm}

\newcommand{\be}{\begin{equation}}
\newcommand{\ee}{\end{equation}}

\newcommand{\klk}{\ifm {,\ldots,} {$,\ldots,$}}

\newcommand{\kpk}{, \ldots ,}

\newenvironment{prf}{\vs\noindent \textbf {Proof}\\}{$\mbox{}$\hfill $\Box$ \\}

\def\cqfd{\vbox{\hrule height 5pt width 5pt }\bigskip} 
\def\qed\cqfd
\def\noi{\noindent}
\def\vs{\smallskip}

\theoremstyle{break}
\newtheorem{proposition}{Proposition}
\newtheorem{definition}[proposition]{Definition} 
 
\newtheorem{theorem}[proposition]{Theorem}
\newtheorem{corollary}[proposition]{Corollary} 
\newtheorem{lemma}[proposition]{Lemma}
 
\newtheorem{observation}[proposition]{Observation}

\newcommand{\ul}{\underline}
\newcommand{\ol}{\overline}
\newcommand{\difrac}{\displaystyle\frac}

\begin{document}
\selectlanguage{english}
\title{\bf {On the geometry of polar varieties}$^{1}$}

\author{\sc B. Bank $^{2}$, M. Giusti $^{3}$, J. Heintz $^{4}$,\\ 
\sc M. Safey El Din $^{5}$, E. Schost $^{6}$}
                 
\maketitle 
\addtocounter{footnote}{1}\footnotetext{Research partially supported
by the following Argentinian, Canadian, French and Spanish agencies
and grants: UBACYT X-098, UBACYT X-113, PICT--2006--02067, the
Canada Research Chair Program, NSERC, BLAN NT05-4-45732 (projet
GECKO), MTM 2007-62799.}
\addtocounter{footnote}{1}\footnotetext{Humboldt-Universit\"at zu Berlin,
Institut f\"ur Mathematik,
D--10099 Berlin, Germany.
bank@mathematik.hu-berlin.de}
\addtocounter{footnote}{1}\footnotetext{CNRS, \'Ecole Polytechnique, Laboratoire LIX,
F--91228 Palaiseau Cedex, France.
Marc.Giusti@Polytechnique.fr}
\addtocounter{footnote}{1}\footnotetext{Departamento de Computaci\'on,
Universidad de Buenos Aires and CONICET, Ciudad Univ., Pab.I, 1428 Ciudad Aut\'onoma de Buenos Aires,
Argentina, and
Departamento de Matem\'aticas,
Estad\'{\i}stica y Computaci\'on, Facultad de Ciencias, Universidad de
Cantabria, 39071 Santander, Spain.\newline 
joos@dc.uba.ar}
\addtocounter{footnote}{1}\footnotetext{UPMC, Univ Paris 06, 
INRIA, Paris-Rocquencourt, SALSA Project; LIP6; CNRS, UMR 7606, LIP6; UFR
Ing\'eni\'erie 919, LIP6 Passy-Kennedy; Case 169, 4, Place Jussieu,
F-75252 Paris \newline Mohab.Safey@lip6.fr}
\addtocounter{footnote}{1}\footnotetext{Computer Science Department,
Room 415, Middlesex College,\\ University of Western Ontario,
Canada.\newline 
eschost@uwo.ca}

\begin{abstract}
\noindent 
We have developed in the past several algorithms 
with intrinsic complexity bounds for the problem of point finding 
in real algebraic varieties. Our aim here is to give 
a comprehensive presentation of the geometrical tools which are necessary to prove the correctness and complexity estimates of these algorithms.
Our results form also the geometrical main ingredients for the
computational treatment of singular hypersurfaces.\\

\noindent
In particular, we show the non--emptiness of suitable generic dual polar varieties
of (possibly singular) real varieties, show that generic polar varieties may become 
singular at smooth points of the original variety and exhibit  a sufficient criterion when this is not the 
case. Further, we introduce the new concept of 
meagerly generic polar varieties and give a degree estimate for them in terms of the degrees of generic polar varieties.\\  The statements are illustrated by examples and a computer
experiment. 

\end{abstract}

\spar


\section{Preliminaries and results}\label{s:1}

Let $\Q$, $\R$ and $\C$ be the fields of rational, real and complex numbers,
respectively, let 
$X:=(\xon)$ be a vector of indeterminates over $\C$ and let be given a reduced
regular sequence $\Fop$ of polynomials in 
$\Q[X]$  such that the ideal $\iFop$ generated by them is the ideal of definition of 
a closed, $\Q$--definable subvariety  $S$ of the $n$--dimensional
complex affine space $\A^n:=\C^n$. Thus $S$ is a non--empty
equidimensional affine variety of dimension $n-p$,
i.e., each irreducible component of $S$ is of dimension $n-p$. 
Said otherwise, $S$ is  of pure codimension $p$ (in $\A^n$).
\spar

We denote by $S_{reg}$ the locus of {\em regular points} of $S$, i.e.,  
the points of $S$, where the Jacobian 
$J(F_1\kpk F_p) :=
\left[\frac{\partial F_k}{\partial X_l}\right]_
{{1 \le k \le p} \atop {1 \le l \le n}}$ has maximal rank $p$, and by 
$S_{sing}:=S\setminus S_{reg}$ the singular locus of $S$.\spar

Let $\A^n_{\R}:=\R^n$ be the $n$--dimensional real affine space. 
We denote by $S_{\R}:=S\cap \A^n_{\R}$ the real trace of the complex variety $S$. 
Moreover, we denote by $\P^n$ the $n$--dimensional complex projective space and by 
$\P^n_{\R}$ its real counterpart. We shall use also the following notations:
\[
S:=\{F_1=0 \klk F_p=0\}\;\;\text{and}\;\;S_{\R}:=\{F_1=0 \klk F_p=0\}_{\R}.
\]
\spar

We denote the coordinate ring of the affine variety $S$ by $\C[S]$. Thus
$\C[S]$ is a finitely generated, reduced, equidimensional $\C$--algebra
which is a domain when $S$ is irreducible.\\
By $\C(S)$ we denote the total quotient ring of $\C[S]$ (or simply of
$S$) which consists of all rational functions of $S$ whose domain has
non--empty intersection with every irreducible component of $S$. When $S$ is
irreducible, then $\C(S)$ becomes the usual field of rational functions
of $S$.\\
The Chinese Remainder Theorem implies that the $\C$--algebra $\C(S)$ is
isomomorphic to the direct product of the function fields of the
irreducible components of $S$.\spar

All varieties that occur in this paper are defined set--theoretically, and not scheme--theoretically. Thus the affine ones have always reduced coordinate rings and when we formulate an algebraic property of a given variety like normality or
Cohen--Macaulayness we refer always to the (reduced) coordinate ring of the variety.

\spar
Let $1\le i \le n-p$ and let 
$\alpha:=\left[ a_{k,l}\right]_{{1 \le k \le n-p-i+1} \atop {0 \le l \le n}}$
be a complex $((n-p-i+1)\times (n+1))$--matrix and suppose that 
$\alpha_{\ast}:=\left[ a_{k,l}\right]_{{1 \le k \le n-p-i+1} \atop {1 \le l \le n}}$ 
has maximal rank $n-p-i+1$.
In case $(a_{1,0}\klk a_{n-p-i+1,0})=0$ we denote by $\ul{K}(\alpha):=\ul{K}^{n-p-i}(\alpha)$ 
and 
in case $(a_{1,0}\klk a_{n-p-i+1,0})\not=0$ by $\ol{K}(\alpha):=\ol{K}^{n-p-i}(\alpha)$  
the $(n-p-i)$--dimensional linear subvarieties of the projective space
$\P^n$ which for $1\le k \le n-p-i+1$ are spanned by the the points 
$(a_{k,0}:a_{k,1}:\cdots : a_{k,n})$.\spar

We define the classic and the dual $i$th polar varieties of $S$ associated 
with the linear varieties  $\ul{K}(\alpha)$ and
$\ol{K}(\alpha)$ as the closures of the loci of the regular points of 
$S$ where all
$(n-i+1)$--minors of the respective polynomial $((n-i+1)\times n)$matrix
\[\begin{bmatrix}
&J\iFop&\\
a_{1,1} & \cdots  & a_{1,n}\\
\vdots & \vdots & \vdots \\
a_{n-p-i+1,1} & \cdots & a_{n-p-i+1,n}\\
\end{bmatrix}
\]
and 
\[\begin{bmatrix}
&J\iFop&\\
a_{1,1}-a_{1,0}X_1 & \cdots  & a_{1,n}-a_{1,0}X_n\\
\vdots & \vdots & \vdots \\
a_{n-p-i+1,1}-a_{n-p-i+1,0}X_1 & \cdots & a_{n-p-i+1,n}-a_{n-p-i+1,0}X_n\\
\end{bmatrix}
\]
vanish.
We denote these polar varieties by 
\[W_{\ul{K}(\alpha)}(S):= 
W_{\ul{K}^{n-p-i}(\alpha)}(S)\;\;\text{and}\;\; W_{\ol{K}(\alpha)}(S):=
W_{\ol{K}^{n-p-i}(\alpha)}(S),
\] 
respectively. They are of expected pure codimension $i$ in $S$. 
Observe also that the polar varieties $W_{\ul{K}^{n-p-i}(\alpha)}(S)$ and
$W_{\ol{K}^{n-p-i}(\alpha)}(S)$ are determined by the $((n-p-i+1)
\times n)$--matrix $a:= {\alpha}_*$. \spar

If $\alpha$ is a real $((n-p-i+1)\times (n+1))$--matrix, we denote by  
\[W_{\ul{K}(\alpha)}(S_{\R}):= W_{\ul{K}^{n-p-i}(\alpha)}(S_{\R}):=
W_{\ul{K}(\alpha)}(S)\cap \A^n_{\R}\]
and
\[W_{\ol{K}(\alpha)}(S_{\R}):=W_{\ol{K}^{n-p-i}(\alpha)}(S_{\R}):=
W_{\ol{K}(a)}(S)\cap \A^n_{\R}\] 
the real traces of $W_{\ul{K}(\alpha)}(S)$ and $W_{\ol{K}(\alpha)}(S)$.
\spar

In this paper we shall work with this purely calculatory definition of the classic and dual polar varieties of $S$. On the other hand, both notions 
may alternatively be charcterized in terms of intrinsic (i.e., coordinate-free), geometric concepts. For example, the classic polar variety $W_{\ul{K}(\alpha)}(S)$
is the  Zariski closure of all regular points $x$ of $S$ such that the tangent space at $x$ is not transversal to the linear variety $\ul{K}(\alpha)$ and the real polar variety $W_{\ul{K}(\alpha)}(S_{\R})$ may be characterized similarly.\spar

In the same vein, fixing an embedding of the affine variety $S$ into the projective space $\P^n$ and fixing a non--degenerate hyperquadric in $\P^n$, we introduced in  \cite{bank3} and \cite{bank4} the notion
of a {\em generalized polar variety} of $S$ which contains as particular instances the notions of classic and dual polar varieties. 
From the intrinsic, geometric characterization of these generalized polar varieties we derived then the present calculatory definition of classic and dual polar varieties.
For details we refer the reader to \cite{bank3} and \cite{bank4}.\spar

The argumentation of this paper will substantially depend on this calculatory definition. This leads to a notation which at first glance looks
overloaded by diacritic marks as e.g. the underscores  and overbars
we use to distinguish classic from dual polar varieties. In view of the context, the reader may overlook in most cases these diacritic marks, however omitting them would introduce an inadmissible amount of ambiguities in statements and proofs.\spar

For  the rest of this section let us assume that for fixed $1\le i \le n-p$ 
there is given a generic complex $((n-p-i+1)\times n)$--matrix 
$a=[a_{k,l}]_{1\le k \le n-p-i+1\atop{1\le l \le n}}$ (the use of the word ''generic'' will be clarified at the end of this section).\\
For $1\le k \le n-p-i+1$ and $0\le l \le n$ we introduce
the following notations:
\[
\ul{a}_{k,l}:=0\;\;\text{and}\;\;\ol{a}_{k,l}:=1\;\;\text{if}\;\;l=0,\;\;\ul{a}_{k,l}:=
\ol{a}_{k,l}:=a_{k,l}\;\;\text{if}\;\;1\le l \le n,
\]
\[
\ul{a}:=[\ul{a}_{k,l}]_{{1\le k \le n-p-i+1}\atop{0\le l \le
n}}\;\;\text{and}\;\;
\ol{a}:=[\ol{a}_{k,l}]_{{1\le k \le n-p-i+1}\atop{0\le l \le n}}
\]
(thus we have $\ul{a}_*=\ol{a}_*=[a_{k,l}]_{{1\le
k \le n-p-i+1}\atop {1\le l \le n}}$).\spar
The corresponding polar varieties
$W_{\ul{K}(\ul{a})}(S)$ and 
$W_{\ol{K}(\ol{a})}(S)$  will be called {\em generic}. 
In this case we have shown in \cite{bank3} and \cite{bank4} that the polar 
varieties $W_{\ul{K}(\ul{a})}(S)$ and $W_{\ol{K}(\ol{a})}(S)$ are either empty or of pure codimension $i$ in $S$. 
Further, we have shown that the local rings of 
$W_{\ul{K}(\ul{a})}(S) $ and $W_{\ol{K}(\ol{a})}(S)$ at any 
regular point are Cohen--Macaulay (\cite{bank3}, Theorem 9).
\spar 

In the case of the classic polar variety $W_{\ul{K}(\ul{a})}(S)$ 
these results are well known to specialists as a consequence of 
Kleiman's variant of the Bertini Theorems (\cite{klei}) and general properties of determinantal varieties (see \cite{p}, 
proof of Lemma 1.3, or \cite{tei1}, Chapitre IV,  proof of Proposition 2). In the case of the dual polar variety 
$W_{\ol{K}(\ol{a})}(S)$ this kind of general argumentation 
cannot be applied and more specific tools are necessary like those developed 
in \cite{bank3} and \cite{bank4}.
\spar

The modern concept of (classic) polar varieties was introduced in the 1930's by 
F. Severi (\cite {se}, \cite{se1}) and J. A. Todd (\cite {to}, \cite{to1}), 
while the intimately related notion of a reciprocal curve goes 
back to the work of J.-V. Poncelet in the period of 1813--1829. 
\spar
As pointed out by Severi and Todd, generic polar varieties have to be understood 
as being organized in certain equivalence classes which embody relevant 
geometric properties of the underlying algebraic variety $S$. 
This view led to the consideration of rational equivalence classes of the 
generic polar varieties. \spar
Around 1975 a renewal of the theory of polar varieties took place with 
essential contributions due R. Piene (\cite {p})  (global theory), B.
Teissier, D. T. L\^e (\cite {lete}, \cite {tei1}), J. P. Henry and M.
Merle (\cite {heme}),  A. Dubson (\cite {du}, Chapitre IV) (local theory),  J. P. Brasselet
and others
(the list is not exhaustive, see \cite{tei2},\cite {p} and \cite {bra} for
a historical account and references).
The idea was to use rational equivalence classes of generic polar varieties 
as a tool which allows to establish numerical formulas in order to 
classify singular varieties by their intrinsic geometric character (\cite {p}).
\spar

 At the same time, classes of generic polar varieties were used in order to formulate 
a manageable local equisingularity criterion which implies the Whitney conditions 
in analytic varieties in view of an intended concept of canonical stratifications 
(see \cite {tei1}).\spar

On the other hand, classic polar varieties became about ten years ago a 
fundamental tool for the design of efficient computer procedures with \
{\em intrinsic} complexity which find real algebraic sample points for 
the connected components of $S_{\R}$, if $\iFop$ is a reduced regular sequence in 
$\Q[X]$ and $S_{\R}$ is smooth and compact.  
The sequential time complexity of these procedures, in its essence, 
turns out to be polynomially bounded by the maximal degree of 
the generic polar varieties $W_{\ul{K}^{n-p-i}(\ul{a})}(S)$, for $1\le i \le n-p$. As we shall see in Section \ref{s:4} as a consequence of Theorem \ref{th:c} , this maximal degree is an invariant of the input system $\iFop$ and the variety $S$, and even of the real variety $S_{\R}$, if all irreducible components of $S$ contain a regular real point. In this sense we refer to the complexity of these algorithms as being intrinsic (see \cite{bank1, bank2}, \cite{mohab2} and \cite{SaSch}) .
\spar

The compactness assumption on $S_{\R}$ was essential in order to 
guarantee the non--emptiness of the classic polar varieties 
$W_{\ul{K}^{n-p-i}(\ul{a})}(S_{\R})$, for $1\le i \le n-p$ . 
If $S_{\R}$ is singular or unbounded, the generic classic polar varieties 
$W_{\ul{K}^{n-p-i}(\ul{a})}(S_{\R})$ may become empty 
(this becomes also a drawback for the geometric analysis of singular varieties 
described above).
\spar

In order to overcome this difficulty at least in the case of a non--singular 
real variety $S_{\R}$, the notion of dual polar varieties was introduced in 
\cite{bank3} and \cite{bank4}.  
The usefulness of dual polar varieties is highlighted by the 
following statement:
\spar

If $S_{\R}$ is smooth, then the dual polar variety 
$W_{\ol{K}^{n-p-i}(\ol{a})}(S_{\R})$ contains a sample point for each 
connected component of $S_{\R}$ (see \cite{bank3, bank4}, Proposition 2 and 
\cite{mopie}, Proposition 2.2).
\spar

In case that $S_{\R}$ is singular, we have the following, considerably 
weaker result which will be shown in Section \ref{s:2} as Theorem \ref{th:a}.
\spar

{\em 
Let $1\le i \le n-p$ and let $C$ be a connected component of the real variety 
$S_{\R}$ containing a regular point. Then, with respect to the 
Euclidean topology, there exists a non--empty, open subset $O^{(i)}_C$ of  
$\A^{(n-p-i+1)\times n}_{\R}$ such that any $((n-p-i+1)\times n)$--matrix 
$a$ of $O^{(i)}_C$
has maximal rank $n-p-i+1$ and such that the real dual polar variety 
$W_{\ol{K}(\ol{a})}(S_{\R})$
is generic and contains a regular point of $C$.
}
\spar

In view of the so--called "lip of Thom" (\cite{thom}), a well studied example 
of a singular curve, this result cannot be improved. Although it is not too expensive 
to construct algorithmically non--empty open conditions which imply the 
conclusion of Theorem \ref{th:a}, the search for rational sample points satisfying 
these conditions seems to be as difficult as the task of finding smooth points 
on singular real varieties.
\spar

Dual polar varieties represent a complex counterpart of the Lagrange multipliers.
Therefore their geometric meaning concerns more real than complex algebraic varieties.
Maybe this is the reason why, motivated by the search for {\em real} solutions           
of polynomial equation systems, they were only recently introduced in (complex) 
algebraic geometry. In the special case of $p:=1$ and $i:=n-p$ the notion of 
a dual polar variety appears implicitly in \cite {Sa05} 
(see also \cite {RoRoSa}, \cite{aumohab} and \cite {Sa07}).\\
The consideration of general $(n-p)$th classic (or dual)
polar varieties was introduced in complexity theory by \cite {grivo} 
and got the name "critical point method". The emerging of elimination
procedures of intrinsic complexity made it necessary to take into account
also the higher dimensional $i$th polar varieties of $S$ (for $1 \le i<n-p$).
\spar

Under the name of "reciprocal polar varieties" generic dual varieties of real 
singular plane curves are exhaustively studied in \cite {mopie} and a manageable 
sufficient condition for their non--emptiness is exhibited.
\spar
An alternative procedure of intrinsic complexity to find sample
points in not necessarily compact smooth semialgebraic varieties was
exhibited in \cite {mohab2}. This procedure is based on the recursive use
of {\em classic} polar varieties.
\spar

We are now going to describe the further content of this paper. For the
sake of simplicity of exposition let us suppose for the moment that  the
given variety $S$ is smooth. Refining
the tools developed in \cite{bank3} and \cite{bank4} we shall show in the
first part of Section \ref{s:3} that the generic
classic and dual polar varieties $W_{\ul{K}(\ul{a})}(S)$ and
$W_{\ol{K}(\ol{a})}(S)$ are normal (see Theorem \ref{th:b} and
Corollary \ref{c:b}). Hence the generic polar varieties of $S$ are both, normal and
Cohen--Macaulay.

Unfortunately, this is the best result we can hope for. In the second part of
Section \ref{s:3} we shall exhibit a general method which allows to obtain smooth varieties $S$ whose  {\em higher} dimensional generic polar varieties are singular.  Hence,  the assertion Theorem 10 (i) of \cite {bank4}, which claims that {\em all} generic polar varieties of $S$ are empty or smooth, is wrong in fact.\\
On the other hand we shall describe a sufficient combinatorial condition in terms of
the parameters $n$, $p$ and $1\le i \le n-p$, that guarantees that the
{\em lower} dimensional generic polar varieties of $S$ are empty or non--singular.
\spar

Let us mention here that in case $p:=1$, i.e., if $S$ is a nonsingular hypersurface,
the classic polar variety
$W_{\ul{K}(\ul{a})}(S)$ is smooth. This is an immediate consequence
of the transversality version of Kleiman's theorem (see also \cite {bank1} for an 
elementary proof). However, in case $p:=1$, the higher dimensional
generic {\em dual} polar varieties of the smooth hypersurface $S$ may
contain singularities.\spar

Finally we explain in Section \ref{s:3} in a more systematic way how singularities
in higher dimensional generic polar varieties of $S$ may arise.\\ Using for
generic $a \in \A^{(n-p-i+1) \times n}$ a natural desingularization of the
(open) polar variety $W_{\ul{K}^{n-p-i}(\ul{a})}(S)\cap S_{reg}$
in the spirit of Room--Kempf  \cite {room, kem}, we shall include the
singularities of $W_{\ul{K}^{n-p-i}(\ul{a})}(S)\cap S_{reg}$ in
a kind of algebraic geometric "discriminant locus" of this
desingularization.
\spar

On the other hand the generic complex $((n-p-i+1)\times n)$--matrix 
$a$ induces an analytic map from $S_{reg}$ to $\A^{n-p-i+1}$. We shall show
that $W_{\ul{K}^{n-p-i}(\ul{a})}(S)\cap S_{reg}$ may be decomposed
into smooth Thom--Boardman strata of this map.\spar

In the geometric analysis of singular varieties as well as in real 
polynomial equation solving,  generic polar varieties play a fundamental role 
as providers of geometric invariants which characterize the underlying 
(complex or real) algebraic variety, in our case $S$ or $S_{\R}$. 
We have already seen that for $1\le i \le n-p$ the generic polar varieties 
$W_{\ul{K}^{n-p-i}(\ul{a})}(S)$ and $W_{\ol{K}^{n-p-i}(\ol{a})}(S)$ 
are empty or of codimension $i$ in $S$ and therefore their dimension becomes 
an invariant of the algebraic variety $S$. On the other hand, Theorem \ref{th:a} 
and the example of Thom's lip imply that for $a\in\Q^{(n-p-i+1)\times n}$ generic 
the dimension of the real polar variety $W_{\ol{K}(\ol{a})}(S_{\R})$ is not an invariant of $S_{\R}$, since $W_{\ol{K}(\ol{a})}(S_{\R})$ may be empty or not, according to the choice of $a$.
\spar

It may occur that the degrees of generic polar varieties represent 
a too coarse measure for the complexity of elimination procedures which solve 
real polynomial equations. Therefore it is sometimes convenient to replace 
for $1\le i \le n-p$ the generic polar varieties of $S$ and $S_{\R}$ by 
more special ones of the form $W_{\ul{K}^{n-p-i}(\ul{a})}(S)$, 
$W_{\ul{K}^{n-p-i}(\ul{a})}(S_{\R})$ and $W_{\ol{K}^{n-p-i}(\ol{a})}(S)$, 
$W_{\ol{K}^{n-p-i}(\ol{a})}(S_{\R})$ (or suitable non--empty Zariski open subsets of them). Here $a$ ranges over a 
Zariski dense  subset of a suitable irreducible subvariety of $\A^{(n-p-i+1)\times n}$ 
(containing generally a Zariski dense set of rational points).
We call these special polar varieties {\em meagerly generic}. 
They share important properties with generic polar varieties 
(e.g. dimension and reducedness) and often they are locally given as transversal intersections of closed form equations and therefore smooth. A particular class of meagerly generic polar varieties with this property was studied in \cite {bank2}.\spar

Another class of meagerly generic polar varieties appears implicitly
in \cite{bank5}, where the problem of finding smooth algebraic sample
points for the (non--degenerated) connected components of {\em singular}
real hypersurfaces is studied.
\spar 

It is not hard to see that the (geometric) degrees of the {\em generic} polar varieties of $S$
constitute invariants of $S$, i.e., the degrees of  $W_{\ul{K}^{n-p-i}(\ul{a})}(S)$
and $W_{\ol{K}^{n-p-i}(\ol{a})}(S)$ are independent of the choice of the (generic) complex $((n-p-i+1)\times n)$--matrix $a$. \spar

The main result of Section \ref{s:4} may be paraphrased as follows: for $1\le i \le n-p$
the degrees of the $i$th meagerly generic classic and dual polar varieties of $S$ are bounded by the degree of the corresponding $i$th generic polar variety of $S$. \spar

The rest of Section \ref{s:4} is devoted to the discussion of the notion of a 
meagerly generic polar variety.\spar

Before finishing this introductory presentation, we add a clarification about our use of the word {\em generic}. To this aim we adopt the point of view of Ren\'e Thom \cite{thom}.
\begin{definition}\label{d:-1}
Let $n,\; 1\le p \le n$ and $1\le i \le n-p$ be as above. The $i$th polar varieties of $S$ depend on linear subvarieties of $\P^n$ which are given by full--rank matrices of the form $a=[a_{k,l}]_{1\le k \le n-p+1\atop{1\le l \le n}}$
with complex, real or rational entries. We say that a given statement is valid for the  generic $i$th classic or dual polar varieties of $S$ if there exists a non--empty
Zariski open (and hence residual dense)  subset $O$ of full--rank matrices of $\A^{(n-p-i+1)\times n}$ (or $\A^{(n-p-i+1)\times n}_{\R}$) such that for any $a$ in $O$ the statement is verified by $W_{\ul{K}(\ul{a})}$ or $W_{\ol{K}(\ol{a})}$ (or their real traces).
\end{definition}

In Section \ref{s:4} we shall face a more general situation: let be given an irreducible
affine subvariety $E$ of $\A^{(n-p-i+1)\times n}$, e.g. an affine linear subspace of
$\A^{(n-p-i+1)\times n}$. With reference to $E$ we shall say that a given statement
is valid for {\em meagerly generic} $i$th polar varieties of $S$ if there exists
a non--empty Zariski open subset $O$ of full--rank matrices of $E$ such that for any $a\in O$ the statement is verified by $W_{\ul{K}(\ul{a})}$ or $W_{\ol{K}(\ol{a})}$(observe that $O$ is residual dense in $E$).\spar 

The aim of this paper is a comprehensive presentation of the geometrical tools which are necessary to prove the correctness of algorithms with intrinsic complexity bounds for real root finding. We developed these algorithms in the past and we think  to develop them further. This leads us to frequent references to already published work on applications of geometric reasoning  to computer science. Thus, in part, this paper has also survey character.

\section {Real dual polar varieties}\label{s:2}

This section is concerned with the proof of Theorem \ref{th:a}, which was announced in 
Section \ref{s:1}. We start with the following technical statement.
 
\begin{lemma}\label{l:a}
Let $C$ be a connected component of the real variety $S_{\R}$ containing an 
regular point.
Then, with respect to the Euclidean topology , there exists a non--empty, open subset
$U_{C}$ of $\A_{\R}^n\setminus S_{\R}$ that satisfies the following condition: 
Let $u$ be an arbitrary point of $U_{C}$ and let $x$ be any point of $S_{\R}$ 
that minimizes the Euclidean distance to $u$ with respect to $S_{\R}$. 
Then $x$ is a regular point belonging to $C$.
\end{lemma}
\begin{prf}
For any two points $z_1,z_2 \in \A^n_{\R}$ and any subset $Y$ of $\A^n_{\R}$ 
we denote by $d(z_1,z_2)$
the Euclidean distance between $z_1$ and $z_2$ and by $d(z_1,Y)$ the 
Euclidean distance from $z_1$ to
$Y$, i.e.,  $d(z_1, Y):=\inf \{ d(z_1,y)\,|\,y\in Y\}$. Let $C_1\klk C_s$ be the 
connected components of 
$S_{\R}$ and suppose without loss of generality that $C=C_1$ holds.
By assumption there exists a regular point $z$ of $C$. 
Observe that the distance 
$d(z, S_{sing}\cap \A^n_{\R})$ 
is positive.\spar

Choose now an open ball $B$ of $\A^n_{\R}$ around the origin which intersects 
$C_1 \klk C_s$
and contains the point $z$ in its interior. 
Since the non--empty sets $C_1\cap \bar{B}\klk C_s\cap \bar{B}$
are disjoint and compact, they have well--defined, positive distances. 
Therefore, there exists a positive real number $d$ stricly smaller than all 
these distances and $d(z, S_{sing}\cap \A^n_{\R})$ such that
the open sphere of radius $d$ and center $z$ is contained in $B$. Let 
$U^*:=\{u^*\in A^n_{\R}\,|\,d(u^*,z)< \frac{d}{2}\}$ be the open ball of 
radius $\frac{d}{2}$ and center $z$ and let $U_C:=U^*\setminus S_{\R}$.\spar

Since $S$ has positive codimension in $\A^n$, we conclude that $U_C$ is a non--empty 
open set with respect to the Euclidean topology which is contained in 
$\A^n_{\R}\setminus S_{\R}$.\spar

Let $u$ be an arbitrary point of $U_C$. Since $u$ does not belong to the closed 
set $S_{\R}$ we have $0< d(u, S_{\R})\le d(u,z)<\frac{d}{2}$. 
Let $x$ be any point of $S_{\R}$ that minimizes the distance to $u$, 
thus satisfying the condition $d(u,x)=d(u,S_{\R})$. From $d(u,x)<\frac{d}{2}$, 
$d(z,u)<\frac{d}{2}$
and the triangle inequality one concludes now that
\[d(z,x)\le d(z,u)+d(u,x)<d\]
holds. Therefore $x$ cannot be contained in $C_2\klk C_s$ and neither in 
$S_{sing}\cap \A^n_{\R}$.
Thus, necessarily $x$ is a regular point of $C$. 
\end{prf}

Now we are going to formulate and prove the main result 
of this section.
\spar

\begin{theorem} \label{th:a}
Let $1\le i \le n-p$ and let $C$ be a connected component of the real variety 
$S_{\R}$
containing a regular point. Then, with respect to the Euclidean topology, 
there exists a non--empty, open subset $O^{(i)}_C$ of  
$\A^{(n-p-i+1)\times n}_{\R}$ such that any $((n-p-i+1)\times n)$--matrix 
$a$ of $O^{(i)}_C$
has maximal rank $n-p-i+1$ and such that the real dual polar variety 
$W_{\ol{K}(\ol{a})}(S_{\R})$
is generic and contains a regular point of $C$.
\end{theorem}
\begin{prf}
Taking into account Lemma \ref{l:a}, we follow the arguments contained in 
the proof of \cite{bank3} and \cite{bank4}, Proposition 2.\spar

Let $C$ be a connected component of $S_{\R}$ containing a regular point and
let $U_C$
be the non--empty open subset of $\A^n \setminus S_{\R}$ introduced in Lemma \ref{l:a}.
Without loss of generality we may assume  $0\notin U_C\setminus S_{\R}$ and that for 
any  point $a\in U_C \setminus S_{\R}$ the dual polar variety 
$W_{\ol{K}^0(a)}(S)$ is empty or generic. 
Thus, putting $O^{(n-p)}_C:=U_C \setminus S_{\R}$, we show first that the statement
of Theorem \ref{th:a} is true for $C$ and $i:=n-p$.\spar

Let $(a_{1,1}\klk a_{1,n})$ be an arbitrary element of $O^{(n-p)}_C$, $a_{1,0}:=1$, 
$a:=(a_{1,0}\klk a_{1,n})$ and let $T^{(n-p)}_a$ be the polynomial 
$((p+1)\times n)$--matrix
\[T^{(n-p)}_a:=T^{(n-p)}_a(X):=\begin{bmatrix}
&J\iFop&\\
a_{1,1}-X_1 & \cdots  & a_{1,n}-X_n\\
\end{bmatrix}.
\]
For $p+1\le k \le n$ denote by $N_k:=N_k(x)$  the $(p+1)$--minor of $T^{(n-p)}_a$
given by the columns $1\klk p,k$. There exists  a point $x$ of $S_{\R}$ that 
minimizes the distance 
to $(a_{1,1}\klk a_{1,n})$ with respect to $S_{\R}$. 
From Lemma \ref{l:a} we deduce that $x$ is
a regular point belonging to $C$. 
Without loss of generality we may assume that
$\det \!\left[\difrac{\partial F_j}{\partial X_k}\right]_{1\le j,k \le p}$ 
does not vanish at $x$.\spar

There exists therefore a chart $Y$ of the real differentiable manifold 
$S_{reg}\cap \A^n_{\R}$,
passing through $x$, with local parameters $X_{p+1|_{Y}}\klk X_{n|_{Y}}$.
Consider now the restriction $\varphi$ of the polynomial function 
$(a_{1,1}-X_1)^2+ \cdots + (a_{1,n}-X_n)^2$ to the chart $Y$ and observe 
that $x\in Y$
minimizes the function $\varphi$ with respect to $Y$. 
From the Lagrange--Multiplier--Theorem \cite{spi} we deduce easily that the \
polynomials $N_{p+1}\klk N_n$
vanish at $x$. Taking into account that 
$\det \!\left[\difrac{\partial F_j}{\partial X_k}\right]_{1\le j,k \le p}$ 
does not vanish at $x$ we infer from the Exchange--Lemma of \cite{bank1} 
that any
$(p+1)$--minor of  $T^{(n-p)}_a$ must vanish at $x$. 
Therefore, $x$ is a regular point belonging to 
$W_{\ol{K}^0(a)}(S_{\R})\cap C$ and $W_{\ol{K}^0(a)}(S)$
is generic.
\spar

Now let $1\le i \le n-p$ be arbitrary and let 
$[a_{k,l}]_{1\le k \le n-p-i+1\atop{1\le l \le n}}$ be a real 
$((n-p-i+1)\times n)$--matrix of maximal rank $n-p-i+1$, 
$a_{1,0}:=\cdots := a_{n-p-i+1,0}:=1$
and $a:=[a_{k,l}]_{1\le k \le n-p-i+1\atop{0\le l \le n}}$.

Further, let $T^{(i)}_a$ be the polynomial $((n-i+1)\times n)$--matrix 
\[T^{(i)}_a:=T^{(i)}_a(X):= \begin{bmatrix}
&J\iFop&\\
a_{1,1}-X_1 & \cdots  & a_{1,n}-X_n\\
\vdots & \vdots & \vdots \\
a_{n-p-i+1,1}-X_1 & \cdots & a_{n-p-i+1,n}-X_n\\
\end{bmatrix}
\]
and recall that the $i$th (complex) dual polar variety of $S$ associated with 
the linear variety
$\ol{K}(a):=\ol{K}^{n-p-i}(a)$, namely $W_{\ol{K}(a)}(S)$, is defined as the 
closure of the locus of the regular points of $S$, where all $(n-p-i)$--minors of 
$T^{(i)}_a$  vanish.\spar  

We choose now a non--empty open subset 
$O^{(i)}_C$ of $\A_{\R}^{(n-p-i+1)\times n}$ satisfying the following conditions:
\begin{itemize}
\item[(i)] Any $((n-p-i+1)\times n)$--matrix  
$[a_{k,l}]_{1\le k \le n-p-i+1\atop{1\le l \le n}}$ belonging to $O^{(i)}_C$
has maximal rank $n-p-i+1$, and for $a_{1,0}:=\cdots := a_{n-p-i+1,0}:=1$
and $a:=[a_{k,l}]_{1\le k \le n-p-i+1\atop{0\le l \le n}}$
the dual polar variety $W_{\ol{K}^{n-p-i}(a)}(S)$
is empty or generic.
\item[(ii)] The point $(a_{1,1}\klk a_{1,n})$  belongs to 
$O^{(n-p)}_C$.
\end{itemize}
Let be given  a $((n-p-i+1)\times (n+1))$--matrix $a$ satisfying 
both conditions above. 
From the structure of the polynomial matrix $T^{(i)}_a$ one infers 
easily that $W_{\ol{K}^0(a)}(S)$ is contained in 
$W_{\ol{K}^{n-p-i}(a)}(S)$. As we have seen above, 
$W_{\ol{K}^0(a)}(S)\cap C$ contains a regular point and 
therefore so does $W_{\ol{K}^{n-p-i}(a)}(S)\cap C$. Moreover, by 
the choice of $O^{(n-p)}_C$,
the dual polar variety $W_{\ol{K}^{n-p-i}(a)}(S)$ is generic.
\end{prf}
\begin{corollary}\label{c:a}
Suppose that the real variety $S_{\R}$ contains a regular point and 
let $1\le i \le n-p$.
Then, with respect to the Euclidean distance, there exists a non-empty,
open subset $O^{(i)}$ of  $\A^{(n-p-i+1)\times n}_{\R}$ such that any 
$((n-p-i+1)\times n)$--matrix 
$[a_{k,l}]_{1\le k \le n-p-i+1\atop{1\le l \le n}}$ of  $O^{(i)}$ has 
maximal rank $n-p-i+1$ and 
such that for  $a_{1,0}:=\cdots := a_{n-p-i+1,0}:=1$
and $a:=[a_{k,l}]_{1\le k \le n-p-i+1\atop{0\le l \le n}}$
the dual polar variety 
$W_{\ol{K}^{n-p-i}(a)}(S)$ is generic and non--empty.  
\end{corollary}
\spar

\section{On the smoothness of generic polar\\ varieties}\label{s:3}

Using the ideas and tools developed in \cite {bank3} and \cite {bank4}, we
are going to prove in the first and main part of this section
that the generic classic and dual polar
varieties of $S$ are normal and Cohen--Macaulay at any point,
where $J\iFop$ has maximal rank. Then we show by means of an infinite
family of examples that it is not always true that {\em all} generic polar
varieties of $S$ are smooth at any point where $J\iFop$ has maximal rank.
\spar 

 We finish the section with two explanations of this phenomenon of non--smoothness.\spar

For the sake of simplicity of exposition, let us assume for the moment
$1\le p < n$ and that any point of $S$ is regular. Further, let
$\left[A_{k,l}\right]_{{1 \le k \le n-p} \atop {1 \le l \le n}}$ be a
$((n-p)\times n)$--matrix of new indeterminates $ A_{k,l}$.
\spar
Recursively in ${1 \le i \le n-p}$ and relative to $\Fop$, we are now
going to introduce two genericity conditions for complex $((n-p)\times
n)$--matrices, namely $U^{(i)}_{classic}$ and $U^{(i)}_{dual}$, such that
\[U_{classic}:= U^{(1)}_{classic}\subset \cdots \subset U^{(n-p)}_{classic}\subset \A^{(n-p)\times n}\]
 and  
\[U_{dual}:= U^{(1)}_{dual}\subset \cdots \subset U^{(n-p)}_{dual}\subset \A^{(n-p)\times n}\]
form two filtrations of $\A^{(n-p)\times n}$ by suitable constructible,
Zariski dense subsets (in fact, we shall choose them as being non--empty and Zariski open).
The sets $U_{classic}$ and $U_{dual}$ will then give an appropriate
meaning to the concept of a generic decreasing sequence  of classic and dual polar
varieties of $S$ or simply to the concept of a generic polar variety. The properties we are going to ensure with this concept of genericity refer to dimension and Cohen--Macaulayness and are used for example in the proofs of Lemma \ref{l:b} and Theorem \ref{th:b} (see Definition \ref{d:-1} for the interpretation of the word "generic" in this paper).
\spar

Let us begin by introducing for ${1 \le i \le n-p}$ the genericity conditions
$U^{(i)}_{classic}$. We start with $i:= n-p$.
\spar

We fix temporarily a $({p \times p})$--minor 
$m$ of the Jacobian $J(F_1\klk F_p)$. For the sake of
conciseness we shall assume
$m := \det \left[\frac{\partial F_j}{\partial X_k}\right]_{{1 \le j \le p} 
\atop {1 \le k \le p}}$. We suppose first that $1\le p < n-1$ holds.
In this case we fix for the moment an arbitrary selection
of indices $1\le k_1\le \cdots \le k_{n-p+1}\le n-p$ and 
$p<l_1\le \cdots \le l_{n-p+1}\le n$,
such that $(k_1,l_1) \klk (k_{n-p+1},l_{n-p+1})$ are all distinct
(observe that the condition $1\le p < n-1$ ensures that such selections
exist). For ${1 \le j \le n-p+1},$ let 
\[ 
M_{k_j,l_j}^{(n-p)} :=\begin{bmatrix}
\frac{\partial F_1}{\partial X_1}&\cdots&\frac{\partial F_1}{\partial X_p}
&\frac{\partial F_1}{\partial X_{l_j}}\\
\vdots&\cdots&\vdots&\vdots\\
\frac{\partial F_p}{\partial X_1}&\cdots&\frac{\partial F_p}{\partial X_p}
&\frac{\partial F_p}{\partial X_{l_j}}\\
A_{k_j,1}&\cdots&A_{k_j,p}&A_{k_j,l_j}\\
\end{bmatrix}.
\]
Writing ${\A^n}_m :=\{x\in \A^n\;|\;m(x)\neq 0\}$,
we may suppose without loss of generality that $S_m := S\cap {\A^n}_m$ is non--empty.\spar
We consider now two polynomial maps of smooth varieties, namely
\[
\Phi:{\A^n}_m \times \A^{(n-p) \times n}\to
\A^n\;\;\text{and}\;\;\Psi:{\A^n}_m\times \A^{(n-p) \times n}\to \A^{n+1}
\]
which are defined as follows:\spar

For $x\in \A^n_m $ and any complex
$((n-p) \times n)$--matrix $a'=[a_{k,l}]_{1\le k \le n-p \atop{1\le l \le n}}$, that contains for each ${1 \le j \le n-p+1}$ at the slots indicated by the
indices $(k_j,1)\klk (k_j,p),$ $(k_j,l_j)$ the entries of the point 
$a_j :=(a_{k_j,1}\klk a_{k_j,p}, a_{k_j,l_j})$ of  $\A^{p+1}$,
the maps $\Phi$ and $\Psi$ take the values
\begin{multline*}
\Phi(x,a'):= \Phi(x,a_1\klk a_{n-p}):=\\(F_1(x)\klk
F_p(x),\,M_{k_1,l_1}(x,a_1)\klk M_{k_{n-p},l_{n-p}}
(x,a_{n-p}))
\end{multline*} and
\begin{multline*}
\Psi(x,a'):= \Psi(x,a_1\klk a_{n-p+1}):=\\(F_1(x)\klk
F_p(x),M_{k_1,l_1}(x,a_1)\klk M_{k_{n-p+1},l_{n-p+1}}(x,a_{n-p+1})).
\end{multline*}
Let $(x,a')$ be an arbitrary point of ${\A^n}_m\times \A^{(n-p) \times n}$
which satisfies the condition $\Phi(x,a')=0$.
Then the Jacobian of $\Phi$ at $(x,a')$ contains a complex $(n \times
(2n-p))$--matrix of the following form: 
\[
\scriptsize\begin{bmatrix}
J\iFop(x)&0&0&\cdots&0\\
              &m(x)&0&\cdots&0\\
               &0&m(x)&\cdots&0\\
        \ast       &\vdots&\vdots&\ddots&\vdots\\
               &&&&\\
               &0&0&\cdots&m(x)
\end{bmatrix}.
\]

Taking into account $m(x)\neq 0$ and that $J(F_1\klk F_p)(x)$ has rank $p$, one sees easily that this matrix has maximal rank $n$. Therefore $\Phi$ is regular
at $(x,a')$. Since this point was chosen arbitrarily in $\Phi^{-1}(0)$,
we conclude that $0\in \A^n$ is a regular value of the polynomial  map $\Phi$.\spar
By a similar argument one infers that $0\in \A^{n+1}$  is a regular value of $\Psi$.\spar
Therefore there exists by the Weak Transversality Theorem of Thom--Sard
(see e.g. \cite{dem}, Ch.3, Theorem 3.7.4, p.79) a non--empty Zariski open subset $U$ of $\A^{(n-p) \times n}$ satisfiying the following conditions:
\spar
For any complex $((n-p)\times n)$--matrix  $a' \in U$ such that $a'$
contains at the slots indicated by the indices $(k_j,1)\klk
(k_j,p),(k_j,l_j)$, ${1 \le j \le n-p+1},$ the entries of suitable points
$a_1\klk a_{n-p+1}$ of $\A^{p+1}$, the equations
\begin{multline}\label{e:c1.1}
F_1(X)=\cdots =F_p(X)=0,\\M_{k_1,\,l_1}(X,a_1)=\cdots =M_{k_{n-p},\,l_{n-p}}(X,a_{n-p})=0
\end{multline} 
intersect transversally at any of their common solutions in $\A^n_m$ and the equations
\begin{multline}\label{e:c1.2}
F_1(X)=\cdots =F_p(X)=0,\\M_{k_1,\,l_1}(X,a_1)=\cdots =M_{k_{n-p+1},\,l_{n-p+1}}(X,a_{n-p+1})=0
\end{multline}
have no common zero in $\A^n_m$.\spar
Remember now that the construction of $U$ depends on the selection of the minor
$m$ and the indices $(k_1,l_1) \klk (k_{n-p+1},l_{n-p+1})$.
There are only finitely many of these choices and each of them gives rise
to a non--empty Zariski open subset of $\A^{(n-p)\times n}$. Cutting the intersection of
all these sets with $(\A^1\setminus \{0\})^{(n-p)\times n}$
we obtain finally $U^{(n-p)}_{classic}$.
\spar
The remaining case $p := n-1$ is treated similarly, considering only the
polynomial map $\Phi$.
\spar
The general step of our recursive construction of the genericity
conditions $U^{(i)}_{classic}$ $1 \le i \le n-p$, is based on the same
kind of argumentation. For the sake of completeness and though
our reasoning may appear repetitive, we shall indicate all essential
points that contain modifications with respect to our previous
argumentation.
\spar
Let ${1 \le i < n-p}$ and suppose that the genericity condition
$U^{(i+1)}_{classic}$ is already constructed. Consider the $(n \times n)$--
matrix
\[
N:=\begin{bmatrix}
&J\iFop&\\
A_{1,1}&\cdots&A_{1,n}\\
\vdots&\vdots&\vdots\\
A_{n-p,1}&\cdots&A_{n-p,n}\\
\end{bmatrix}
\] 
and fix for the moment an arbitrary $((n-i) \times (n-i))$-- submatrix of
$N$ which contains $n-i$ entries from each of the rows number $1\klk p$ of
$N$. Let $m$ denote the corresponding $(n-i)$--minor of $N$. For the sake of
conciseness we shall assume that $m$ is the minor
\[
m:=\det \begin{bmatrix}
\frac{\partial F_1}{\partial X_1}&\cdots&\frac{\partial F_1}{\partial X_{n-i}}\\
\vdots&\cdots&\vdots\\
\frac{\partial F_p}{\partial X_1}&\cdots&\frac{\partial F_p}{\partial X_{n-i}}\\
A_{1,1}&\cdots&A_{1,n-i}\\
\vdots&\vdots&\vdots\\
A_{n-p-i,1}&\cdots&A_{n-p-i,n-i}\\
\end{bmatrix}.
\]
Writing $A:=[A_{k,l}]_{1\le k \le n-p-i\atop{1\le l \le n-i }}$ we denote for $a''$
in $\A^{(n-p-i)\times n}$ by $m(X,a'')$ the polynomial obtained from $m=m(X,A)$ by specializing $A$ to the complex $((n-p-i)\times n)$--matrix $a''$.\\
Without loss of generality we may suppose that $(S \times (\A^{(n-p)\times n})_m$ is non--empty. \spar
Let us first suppose that $i^2 > n-p $ holds.\spar
In this case we fix temporarily an arbitrary selection of indices
$n-p-i < k_1\le \cdots \le k_{n-p+1}\le n-p$ and
$n-i < l_1\le \cdots \le l_{n-p+1}\le n$
such that $(k_1,l_1) \klk (k_{n-p+1},l_{n-p+1})$ are all distinct
(observe that the condition $i^2 > n-p$ ensures that such selections exist).
\spar

For $1 \le j \le n-p+1$ we shall consider in the following the
$(n-i+1)$--minor 
\[
M_{k_j,l_j}^{(i)} :=\det \begin{bmatrix}
\frac{\partial F_1}{\partial X_1}&\cdots&\frac{\partial F_1}{\partial X_{n-i}}&\frac{\partial F_1}{\partial X_{l_j}}\\
\vdots&\cdots&\vdots&\vdots\\
\frac{\partial F_p}{\partial X_1}&\cdots&\frac{\partial F_p}{\partial X_{n-i}}
&\frac{\partial F_p}{\partial X_{l_j}}\\
A_{1,1}&\cdots&A_{1,n-i}&A_{1,l_j}\\
\vdots&\vdots&\vdots&\vdots\\
A_{n-p-i,1}&\cdots&A_{n-p-i,n-i}&A_{n-p-i,l_j}\\
A_{k_j,1}&\cdots&A_{k_j,n-i}&A_{k_j,l_j}\\
\end{bmatrix}.
\] 
of $N$.

In the same spirit as before, we introduce now two polynomial maps of smooth varieties, namely
\[
\Phi: {(\A^n \times \A^{(n-p-i)\times n})_m} \times
\A^{i \times n} \to \A^n
\]
and
\[\qquad \Psi: {(\A^n \times
\A^{(n-p-i)\times n})_m} \times \A^{i \times n}\to \A^{n+1}
\] 
which are defined as follows:
\spar

For $(x,a'')\in {(\A^n \times \A^{(n-p-i)\times n})_m}$ and any complex
$(i \times n)$--matrix \\$a''':=[a_{k,l}]_{n-p-i < k \le n-p \atop{1 \le l \le n}}$ 
that contains for ${1 \le j \le n-p+1}$
at the slots indicated by the indices $(k_j,1)\klk (k_j,n-i),(1,l_j) \klk (n-p-i,l_j),(k_j,l_j)$ the entries of the point
$a_j :=(a_{k_j,1}\klk a_{k_j,n-i},a_{1,l_j}\klk a_{n-p-i,l_j},a_{k_j,l_j})$  of 
$ \A^{2(n-i)-p+1}$, the maps $\Phi$ and $\Psi$ take the values 
\begin{multline*}
\Phi(x,a'',a'''):=\Phi(x,a'',a_1\klk a_{n-p}):=\\(F_1(x)\klk
F_p(x),M_{k_1,l_1}(x,a'',a_1)\klk M_{k_{n-p},l_{n-p}}
(x,a'',a_{n-p}))
\end{multline*}
and
\begin{multline*}
\Psi(x,a'',a'''):=\Psi(x,a'',a_1\klk a_{n-p+1}):=\\
:=(F_1(x)\klk F_p(x),M_{k_1,l_1}(x,a'',a_1)\klk
M_{k_{n-p+1},l_{n-p+1}}(x,a'',a_{n-p+1})).
\end{multline*}
For a fixed full rank matrix $a''\in \A^{(n-p-i)\times n}$ we denote by 
\[
\Phi_{a''}:\A^n_{m(X,a'')}\times  \A^{i\times n}\to \A^n\;\;\text{and}\;\;\
\Psi_{a''}:\A^n_{m(X,a'')}\times  \A^{i\times n}\to \A^{n+1}
\]
the polynomial maps of smooth varieties induced by $\Phi$ and $\Psi$ when we specialize the second argument to the value $a''$.
\spar

Let $(x,a'',a''')$ be an arbitrary point of ${(\A^n \times
\A^{(n-p-i)\times n})_m} \times \A^{i \times n}$
which satisfies the condition $\Phi(x,a'',a''')=0$.
Then the Jacobian of $\Phi_{a''}$ at $(x,a''')$ contains a complex 
$(n\times (2n-p))$--matrix of the following form
\[
\scriptsize\begin{bmatrix}
J\iFop(x)&0&0&\cdots&0\\
              &m(x,a'')&0&\cdots&0\\
               &0&m(x,a'')&\cdots&0\\
        \ast       &\vdots&\vdots&\ddots&\vdots\\
               &&&&\\
               &0&0&\cdots&m(x,a'')
\end{bmatrix}.
\]
Taking into account $m(x,a'')\neq 0$ and that $J(F_1\kpk F_p)(x)$
has rank $p$, one sees easily that this matrix and hence the
Jacobian of $\Phi_{a''}$ at the point $(x,a''')$ has maximal rank $n$.
Therefore $\Phi_{a''}$ is regular at $(x,a''')$. Since this point was chosen arbitrarily in
$\Phi_{a''}^{-1}(0)$, we conclude that $0\in \A^n$ is a regular value of the polynomial  map $\Phi_{a''}$. Similarly one argues that $0\in \A^{n+1}$  is a regular value of
$\Psi_{a''}$.\spar
We consider now the constructible subset $U$ of
$\A^{(n-p)\times n}$ defined by the following conditions:\spar

For any pair consisting of a complex $((n-i) \times n)$--matrix $a''$ and a complex $(i \times n)$--matrix $a'''$ such that 
$(a'',a''')$ belongs to $U$ and such that $a'''$
contains at the slots indicated by the indices $(k_j,1)\klk
(k_j,n-i),(1,l_j)\klk (n-p-i,l_j),(k_j,l_j)$, ${1 \le j \le n-p+1}$,
the entries of suitable points
$a_1\klk a_{n-p+1}$ of $\A^{2(n-i)-p+1}$, the equations
\begin{multline}\label{e:c2.1}
F_1(X)=\cdots =F_p(X)=0, \\
M_{k_1,\,l_1}(X,a'',a_1)=\cdots
=M_{k_{n-p},\,l_{n-p}}(X,a'',a_{n-p})=0
\end{multline}
intersect transversally at any of their common solutions in $\A^n_{m(X,a'')}$ and the equations
\begin{multline}\label{e:c2.2}
F_1(X)=\cdots =F_p(X)=0,\\
M_{k_1,\,l_1}(X,a'',a_1)=\cdots
=M_{k_{n-p+1},\,l_{n-p+1}}(X,a'',a_{n-p+1})=0
\end{multline}
have no common zero in $\A^n_{m(X,a'')}$.
\spar
Applying now for every full rank matrix $a''\in \A^{(n-p-i)\times n}$ the Weak Transversality Theorem of Thom--Sard to $\Phi_{a''}$ and $\Psi_{a''}$, we conclude from the constructibility
of $U$ that it is Zariski dense in $\A^{(n-p)\times n}$. Therefore we may suppose 
without loss of generality that $U$ is a non--empty Zariski open subset of $\A^{(n-p)\times n}$ (here and in the sequel we use the fact that a residual dense, constructible subset of an affine space contains a non--empty Zariski open subset).

\spar
Remember now that the construction of $U$ depends on the selection of the
minor
$m$ and the indices $(k_1,l_1) \klk (k_{n-p+1},l_{n-p+1})$.
There are only finitely many of these choices and each of them gives rise
to a non--empty Zariski open subset of $\A^{(n-p)\times n}$. Cutting the
intersection of all these sets with $U^{(i+1)}_{classic}$  we obtain
finally
$U^{(i)}_{classic}$.
\spar
The remaining case  $i^2 \le n-p$ is treated similarly. In order to
explain the differences with the previous argumentation, let us fix again
the $(n-i)$-- minor 
\[
m:=\det \begin{bmatrix}
\frac{\partial F_1}{\partial X_1}&\cdots&\frac{\partial F_1}{\partial X_{n-i}}\\
\vdots&\cdots&\vdots\\
\frac{\partial F_p}{\partial X_1}&\cdots&\frac{\partial F_p}{\partial X_{n-i}}\\
A_{1,1}&\cdots&A_{1,n-i}\\
\vdots&\vdots&\vdots\\
A_{n-p-i,1}&\cdots&A_{n-p-i,n-i}\\
\end{bmatrix}
\] of the
$n\times n$ matrix $N$, and an arbitrary selection of indices $n-p-i <
k_1\le \cdots \le k_{i^2}\le n-p$ and
$n-i < l_1\le \cdots \le l_{i^2}\le n$,
such that $(k_1,l_1) \klk (k_{i^2},l_{i^2})$ are all distinct
(observe that the condition $i^2 \le n-p$ ensures that such selections
exist).\\
The $(n-i+1)$-- minors $M_{k_j,l_j},\;1 \le j \le i^2$ are the same as
before. We consider now instead of $\Phi$ and $\Psi$ only the polynomial
map
\[
\widetilde{\Phi} : (\A^n \times \A^{(n-p-i)\times n})_m \times
\A^{i \times n}\to \A^{p+i^2}
\]
which is defined as follows:
\spar

For $(x,a'')\in (\A^n \times \A^{(n-p-i)\times n})_m$
and any complex $(i \times n)$--matrix \\
$a''':=[a_{k,l}]_{n-p-i< k \le n-p\atop{1\le l \le n}}$ that contains for $1 \le j \le i^2$
at the slots indicated by the
indices $(k_j,1)\klk (k_j,n-i),(1,l_j)\klk (n-p-i,l_j),(k_j,l_j)$ the entries of the point
$a_j :=(a_{k_j,1}\klk a_{k_j,n-i},a_{1,l_j}\klk a_{n-p-i,l_j},a_{k_j,l_j})$
of $\A^{2(n-i)-p+1}$,
the map $\widetilde{\Phi}$ takes the value
\begin{multline*}
\widetilde{\Phi}(x,a'',a'''):=\widetilde{\Phi}(x,a'',a_1\klk a_{i^2}):=\\
(F_1(x)\klk
F_p(x),M_{k_1,l_1}(x,a'',a_1)\klk M_{k_{i^2},l_{i^2}}
(x,a'',a_{i^2})).
\end{multline*}
Similarly as before we denote for a fixed full rank matrix $a''\in \A^{(n-p-i)\times n}$ by
\[
\widetilde{\Phi}_{a''}:\A^n_{m(X,a'')}\times \A^{i\times n}\to \A^{p+i^2}
\]
the polynomial map of smooth varieties induced by $\widetilde{\Phi}$ when we specialize the second argument to the value $a''$. Then we conclude again that for any 
$a''\in \A^{(n-p-i)\times n}$ the point $0\in \A^{p+i^2}$ is a regular value of the
polynomial map $\widetilde{\Phi}_{a''}$ and that there exists a non--empty Zariski
open subset $\widetilde{U}$ of $\A^{(n-p)\times n}$ satisfying the following
condition:\spar
For any pair consisting of a complex $((n-i) \times n)$--matrix $a''$
and a complex $(i \times n)$--matrix $a'''$ such that
$(a'',a''')$ belongs to $\widetilde{U}$ and  such that $a'''$
contains at the slots indicated by the indices $(k_j,1)\klk
(k_j,n-i),(1,l_j)\klk (n-p-i,l_j),(k_j,l_j),\;1 \le j \le i^2$, the
entries of suitable points $a_1 \klk a_{i^2}$ of $\A^{2(n-i)-p+1}$, the equations
\begin{multline}\label{e:c3}
F_1(X)=\cdots =F_p(X)=0, \\
M_{k_1,\,l_1}(X,a'',a_1)=\cdots
=M_{k_{n-p},\,l_{n-p}}(X,a'',a_{i^2})=0
\end{multline}
intersect transversally at any of their common solutions in $\A^n_{m(X,a'')}$.
\spar
Replacing in the previous argumentation the set $U$ by $\widetilde{U}$ we
define now $U^{(i)}_{classic}$ in the same way as before.
\spar
The construction of a filtration of $\A^{(n-p)\times n}$ by non-empty
Zariski open subsets $U^{(i)}_{dual},\;\;1 \le j \le n-p$, follows the same
line of reasoning, where the $n \times n$-- matrix $N$
has to be replaced by 
\[
\begin{bmatrix}
&J\iFop&\\
A_{1,1}-X_1&\cdots&A_{1,n}-X_n\\
\vdots&\vdots&\vdots\\
A_{n-p,1}-X_1&\cdots&A_{n-p,n}-X_n\\
\end{bmatrix}.
\]
One has only to take care to add in the construction of the sets
$U^{(i)}_{dual}$, $1 \le i \le n-p$, suitable Zariski open 
conditions for the minors of the $((n-p)\times n)$-- matrix 
\[
[A_{k,l}]_{1\le k \le n-p \atop{1\le l \le n}}.
\]
For the rest of this section we fix a complex $((n-p)\times n)$--matrix 
$[a_{k,l}]_{{1\le k \le n-p}\atop {1\le l \le n}}$ belonging to the genericity condition
$U_{classic}\cap U_{dual}$.\spar
For $1\le i \le n-p,\;\;1\le k \le n-p-i+1,\;\;0\le l \le n$ we introduce the following notations:
\[
\ul{a}^{(i)}_{k,l}:=0\;\;\text{and}\;\;\ol{a}^{(i)}_{k,l}:=1\;\;\text{if}\;\;l=0,\;\;\ul{a}^{(i)}_{k,l}=
\ol{a}^{(i)}_{k,l}:=a_{k,l}\;\;\text{if}\;\;1\le l \le n,
\]
\[
\ul{a}^{(i)}:=[\ul{a}^{(i)}_{k,l}]_{{1\le k \le n-p-i+1}\atop{0\le l \le n}}\;\;\text{and}\;\;
\ol{a}^{(i)}:=[\ol{a}^{(i)}_{k,l}]_{{1\le k \le n-p-i+1}\atop{0\le l \le n}}
\]
(thus we have in terms of the notation of
Section \ref{s:1} the identity 
${\ul{a}^{(i)}}_*={\ol{a}^{(i)}}_*=[a_{k,l}]_{{1\le k \le n-p-i+1}\atop {1\le l \le n}}).$

Since $[a_{k,l}]_{{1\le k \le n-p-i+1}\atop {1\le l \le n}}$ belongs by assumption to the non--empty Zariski open subset $U_{classic}\cap U_{dual}$ of $\A^{(n-p)\times n}$, we shall consider for $1\le i \le n-p$ the matrices $\ul{a}^{(i)}$ and $\ol{a}^{(i)}$ and the corresponding classic and dual polar varieties 
\[
W_{\ul{K}(\ul{a}^{(i)})}(S)=W_{\ul{K}^{n-p-i}(\ul{a}^{(i)})}(S) \;\;\text{and}\;\; 
W_{\ol{K}(\ol{a}^{(i)})}(S)=W_{\ol{K}^{n-p-i}(\ol{a}^{(i)})}(S)\]
as ''generic''. These polar varieties are organized in two nested sequences
\[
W_{\ul{K}(\ul{a}^{(n-p)})}(S)\subset \cdots \subset W_{\ul{K}(\ul{a}^{(1)})}(S)\subset S \subset \A^n\]
and 
\[
W_{\ol{K}(\ol{a}^{(n-p)})}(S)\subset \cdots \subset W_{\ol{K}(\ol{a}^{(1)})}(S)\subset S \subset \A^n.
\]
Since by construction $U_{classic}$ satisfies the conditions (\ref{e:c1.1}),
(\ref{e:c1.2}), (\ref{e:c2.1}), (\ref{e:c2.2}), (\ref{e:c3}) and $U_{dual}$ behaves
mutatis mutandis in the same way, we conclude that for $1\le i \le n-p$ the polar varieties
$W_{\ul{K}(\ul{a}^{(i)})}(S)$ and $W_{\ol{K}(\ol{a}^{(i)})}(S)$
are also generic in the sense of \cite{bank3} and \cite{bank4}.\spar
Let $1\le i < n-p$ and $1\le h \le n-p-i+1$. We denote by $\ul{a}^{(i,h)}$ and $\ol{a}^{(i,h)}$ the $((n-p-i)\times (n+1))$--matrices obtained by from $\ul{a}^{(i)}$ and
$\ol{a}^{(i)}$ by deleting their row number $h$, namely 
\[\ul{a}^{(i,h)}:=[\ul{a}^{(i)}_{k,l}]_{{1\le k \le n-p-i+1}\atop{{k\neq h}\atop{0\le l \le n}}}\quad \text{and}\quad \ol{a}^{(i,h)}:=[\ol{a}^{(i)}_{k,l}]_{{1\le k \le n-p-i+1}\atop{{k\neq h}\atop{0\le l \le n}}}.
\]
Thus we have, in particular, 
\[
\ul{a}^{(i+1)}=\ul{a}^{(i,n-p-i+1)}\;\;\text{and}\;\; \ol{a}^{(i+1)}=\ol{a}^{(i,n-p-i+1)}.
\]
According to the notations introduced in \cite{bank3} and \cite{bank4} we write
\[
\ul{\Delta}_i:=\bigcap_{1\le h \le n-p-i+1}W_{\ul{K}(\ul{a}^{(i,h)})}(S)
\;\;\text{and}\;\;\ol{\Delta}_i:=\bigcap_{1\le h \le n-p-i+1}W_{\ol{K}(\ol{a}^{(i,h)})}(S).
\] 
From \cite{bank4} Proposition 6 and the subsequent commentaries we conclude that the
singular loci of 
$W_{\ul{K}(\ul{a}^{(i)})}(S)$ and $W_{\ol{K}(\ol{a}^{(i)})}(S)$
are contained in $\ul{\Delta}_i$ and $\ol{\Delta}_i$, respectively.\spar
The crucial point of the construction of the genericity conditions 
 $U_{classic}$ and  $U_{dual}$ may be summarized in the following statement. 
\begin{lemma}\label{l:b}
Let $1\le i < n-p$. Then  $\ul{\Delta}_i$ and $\ol{\Delta}_i$ are empty or closed subvarieties of 
$W_{\ul{K}(\ul{a}^{(i)})}(S)$ and $W_{\ol{K}(\ol{a}^{(i)})}(S)$
of codimension $\ge 2$, respectively.
\end{lemma}
\begin{prf}
For the sake of simplicity of exposition, we restrict our attention to $\ul{\Delta}_i$.
The case of $\ol{\Delta}_i$ can be treated in the same way.\spar

Let us first assume $i=n-p-1$. Suppose that  $\ul{\Delta}_{n-p-1}$ is non--empty and consider an arbitrary point $x\in \ul{\Delta}_{n-p-1}$. Since $S$ is smooth, there is a $p$--minor $m$ of $J\iFop$, say
\[
m:=\det\;\begin{bmatrix}
\frac{\partial F_1}{\partial X_1}&\cdots&\frac{\partial F_1}{\partial X_p}\\
\vdots&\cdots&\vdots\\
\frac{\partial F_p}{\partial X_1}&\cdots&\frac{\partial F_p}{\partial X_p}\\
\end{bmatrix}
\]
with $m(x)\neq0$. Since $x$ belongs to $\ul{\Delta}_{n-p-1}$, we have
\[F_1(x)=0\; \klk  F_p(x)=0\]
and
\[ M_{1,p+1}(x,a^{(n-p-1, 2)})=0\; \klk M_{1,n}(x,a^{(n-p-1, 2)})=0,\;\;
M_{2,p+1}(x,a^{(n-p-1, 1)})=0,\]
in contradiction to the conditions (\ref{e:c1.1}), (\ref{e:c1.2}) satisfied by $U_{classic}$.
Therefore $\ul{\Delta}_{n-p-1}$ must be empty. This proves Lemma \ref{l:b} in case
$i:=n-p-1$.\spar

Now suppose $1\le i <n-p-1$ and that $\ul{\Delta}_{i+1}$ is either empty or of codimension 
$\ge 2$ in $W_{\ul{K}(\ul{a}^{(i+1)})}(S)$. Let
\[
N_i:=\begin{bmatrix}
&J\iFop&\\
a_{1,1}&\cdots&a_{1,n}\\
\vdots&\vdots&\vdots\\
a_{n-p-i,1}&\cdots&a_{n-p-i,n}\\
\end{bmatrix}.
\] 
If $\ul{\Delta}_i$ is empty we are done. Otherwise consider an arbitrary irreducible component $C$ of $\ul{\Delta}_i$. Assume first that for any $((n-i-1)\times (n-i-1))$--
submatrix of $N_i$, which contains $n-i-1$ entries from each of the rows number 
$1\klk p$ of $N_i$, the determinant vanishes identically on $C$. This implies that $C$ is
contained in  $\ul{\Delta}_{i+1}$. From our assumptions on $\ul{\Delta}_{i+1}$ we deduce that $C$ must be of codimension $\ge 2$ in $W_{\ul{K}(\ul{a}^{(i+1)})}(S)$
and hence in $W_{\ul{K}(\ul{a}^{(i)})}(S)$. \spar

Therefore we may suppose without loss of generality that there exists a $((n-i-1)\times (n-i-1))$--submatrix of $N_i$ containing $n-i-1$ entries from each of the rows number $1 \klk p$ of $N_i$, whose determinant, say $m$, does not vanish identically on $C$.\spar

For the sake of conciseness  we shall assume  
\[
m:=\det\;\begin{bmatrix}
\frac{\partial F_1}{\partial X_1}&\cdots&\frac{\partial F_1}{\partial X_{n-i-1}}\\
\vdots&\cdots&\vdots\\
\frac{\partial F_p}{\partial X_1}&\cdots&\frac{\partial F_p}{\partial X_{n-i-1}}\\
a_{1,1}&\cdots&a_{1,n-i-1}\\
\vdots&\cdots&\vdots\\
a_{n-p-i-1,1}&\cdots&a_{n-p-i-1,n-i-1}\\
\end{bmatrix}.
\]
Since $C$ is contained in $\ul{\Delta}_i$, any point of $C$ is a zero of 
the equation system
\begin{multline}\label{e:1}
F_1(X)=0 \klk F_p(X)=0,\\
M_{n-p-i,n-i}(X,a^{(i,n-p-i+1)})=0\; \klk M_{n-p-i,n}(X,a^{(i, n-p-i+1)})=0,\\
M_{n-p-i+1,n-i}(X,a^{(i,n-p-i)})\; \klk M_{n-p-i+1,n}(X,a^{(i,
n-p-i)})=0.
\end{multline}
Since $U_{classic}$ satisfies the conditions (\ref{e:c2.1}), (\ref{e:c2.2}) 
and (\ref{e:c3}),
we conclude that the equations of (\ref{e:1}) intersect transversally
 at any point of the (non--empty) set $C_m$. 
 This implies $\dim C \le n-p-2(i+1)$.
\spar

On the other hand we deduce from \cite{bank3}, \cite{bank4}, 
Proposition 8 that $\dim W_{\ul{K}(\ul{a}^{(i)})}(S)=n-p-i$ holds.  
Since $C$ is a closed irreducible subvariety of 
$W_{\ul{K}(\ul{a}^{(i)})}(S)$,
we infer that the codimension of $C$ in $W_{\ul{K}(\ul{a}^{(i)})}(S)$ is at least $(n-p-i)-(n-p-2(i+1))=i+2\ge2.$
\end{prf}
From Lemma \ref{l:b} we draw now the following conclusion.
\begin{theorem}\label{th:b}
Suppose that $S$ is smooth and let $1\le i \le n-p$. Then the generic 
polar varieties $W_{\ul{K}(\ul{a}^{(i)})}(S)$ and
$W_{\ol{K}(\ol{a}^{(i)})}(S)$ are empty or normal Cohen--Macaulay 
subvarieties of $S$ of pure  codimension $i$.
\end{theorem}
\begin{prf}
We may restrict our attention to $W_{\ul{K}(\ul{a}^{(i)})}(S)$. The case of 
$W_{\ol{K}(\ol{a}^{(i)})}(S)$ can be treated in a similar way.\spar
From  \cite{bank3}, Theorem 9 and Corollary 10 we deduce that $W_{\ul{K}(\ul{a}^{(i)})}(S)$ is empty or a Cohen--Macaulay subvariety of $S$ of codimension $i$.\spar
We are now going to show that $W_{\ul{K}(\ul{a}^{(i)})}(S)$ is normal.\spar

From \cite{bank3}, \cite{bank4}, Lemma 7 we deduce immediately that $W_{\ul{K}(\ul{a}^{(n-p)})}(S)$ is empty or $0$--dimensional and hence normal
(recall from Section \ref{s:1} that $W_{\ul{K}(\ul{a}^{(n-p)})}(S)$ is defined set--theoretically and that its coordinate ring  is therefore reduced).
\spar

Therefore we may assume without loss of generality $1\le i < n-p$ and $W_{\ul{K}(\ul{a}^{(i)})}(S)\neq \emptyset.$ As observed above, the singular locus of
$W_{\ul{K}(\ul{a}^{(i)})}(S)$ is contained in $\ul{\Delta}_i$ and hence, by Lemma \ref{l:b}, empty or of codimension $\ge2$ in $W_{\ul{K}(\ul{a}^{(i)})}(S)$. Consequently, $W_{\ul{K}(\ul{a}^{(i)})}(S)$ is regular in codimension one.\spar
Since $W_{\ul{K}(\ul{a}^{(i)})}(S)$ is Cohen--Macaulay, we infer from
Serre's normality criterion (see e.g. \cite{mat}, Theorem 23.8) that 
$W_{\ul{K}(\ul{a}^{(i)})}(S)$ is a normal variety.
\end{prf}
\begin{observation}\label{o:a}
From  \cite{bank3}, Theorem 9 one deduces that the defining ideal of $W_{\ul{K}(\ul{a}^{(i)})}(S)$ in $\Q[X]$ is generated by $\Fop$ and all $(n-i+1)$--minors of the polynomial matrix $((n-i+1)\times n)$--matrix 
\[
\begin{bmatrix}
&J\iFop&\\
a_{1,1}&\cdots&a_{1,n}\\
\vdots&\vdots&\vdots\\
a_{n-p-i+1,1}&\cdots&a_{n-p-i+1,n}\\
\end{bmatrix}.
\] 
\end{observation}
\spar
For the rest of this section let us drop the assumption that the variety $S$ is smooth. Then Theorem \ref{th:b} may be reformulated as follows.\spar

\begin{corollary}\label{c:b}
Let $1\le i \le n-p$. Then the generic polar varieties $W_{\ul{K}(\ul{a}^{(i)})}(S)$ and $W_{\ol{K}(\ol{a}^{(i)})}(S)$ are empty or subvarieties of $S$ of pure codimension $i$ that are normal and Cohen--Macaulay at any point where $J\iFop$ has maximal rank.

\end{corollary}
\begin{prf}
Let $m$ be an arbitrary $p$--minor of the Jacobian $J\iFop$. The proof of Theorem \ref{th:b} relies only on the assumption that $S$ is smooth. Its correctness does not require that
$S$ is closed. Therefore the statement of Theorem \ref{th:b} remains mutatis mutandis correct if we replace $S$ by $S_m$. Therefore $S_m\neq \emptyset$ implies that
\[(W_{\ul{K}(\ul{a}^{(i)})}(S))_m=W_{\ul{K}(\ul{a}^{(i)})}(S_m)\;\;\text{and}\;\;(W_{\ol{K}(\ol{a}^{(i)})}(S))_m=W_{\ol{K}(\ol{a}^{(i)})}(S_m)\]
are empty or normal Cohen--Macaulay subvarieties of $S_m$ of codimension $i$. Since $m$ was an arbitary $p$--minor of $J\iFop$ the assertion of Corollary \ref{c:b} follows immediately.
\end{prf}

Here the following geometric comment is at order. Schubert varieties are normal and Cohen--Macaulay, whereas classic polar varieties may be interpreted as pre--images of Schubert varieties under suitable polynomial maps (see \cite{p}, proof of Lemma 1.3 (ii)).  This fact suggests the statements of Theorem \ref{th:b} and Corollary \ref{c:b} and an alternative proof of them. In the dual case these statements are new. \spar

We deduced for  $a\in \A^{(n-p)\times n}$ generic and $1\le i \le n-p$ the normality of the polar varieties $W_{\ul{K}(\ul{a}^{(i)})}(S)$ and $W_{\ol{K}(\ol{a}^{(i)})}(S)$ of $S$ from the consideration of the Zariski closed subsets $\ul{\Delta}_i$ and $\ol{\Delta}_i$ of $\A^n$ which contain the singular locus of $W_{\ul{K}(\ul{a}^{(i)})}(S)$ and $W_{\ol{K}(\ol{a}^{(i)})}(S)$.
The main tool was Lemma \ref{l:b} which in its turn relied on the elementary, but rather tedious construction of the non--empty Zariski open subsets $\mathcal{U}_{classic}$ and $\mathcal{U}_{dual}$
of $\A^{(n-p)\times n}$. It is necessary to invest some work in proofs because
$\ul{\Delta}_i$ 
is, unlike the polar varieties, not determinantal.\spar

General statements describing the singular locus of a determinantal variety as another determinantal variety (see e.g. \cite{bru}) cannot be applied in a straightforward manner to polar varieties. The problem arises from the fact that this property of determinantal
varieties is generally not preserved under pre--images. 
However, as it was pointed out by one of the referees, the classic polar varieties can be interpreted as degeneracy loci of vector bundle maps (see Examples 14.3.3 and 14.4.15 in \cite{fu}). Thus the global point
of view of  degeneracy loci of bundle maps may be applied to this situation and it can be shown that the singular locus of the generic classic polar variety $W_{\ul{K}(\ul{a}^{(i)})}(S)$
is empty or of dimension $n-p-(2i+2)$. This implies Theorem \ref{th:b} and the next
proposition in the classic case. In the dual case both statements are undoubtedly new.\spar

We discuss now under which conditions it may be guaranteed that a generic polar variety of $S$ is smooth at any point where the $J\iFop$ has maximal rank.\spar

\begin{proposition}\label{p:a}
Let $1\le i \le n-p$ with $2i+2>n-p$. Then the generic polar varieties $W_{\ul{K}(\ul{a}^{(i)})}(S)$ and $W_{\ol{K}(\ol{a}^{(i)})}(S)$ are smooth at any point
where $J\iFop$ has maximal rank.
\end{proposition}
\begin{prf}
Again we restrict our attention to $W_{\ul{K}(\ul{a}^{(i)})}(S)$. The case of $W_{\ol{K}(\ol{a}^{(i)})}(S)$ is treated similarly. \spar

From \cite{bank4}, Proposition 6 and the subsequent commentaries we conclude that the points at which $W_{\ul{K}(\ul{a}^{(i)})}(S)$ is singular
and $J\iFop$ has maximal rank
are contained in $\ul{\Delta}_i$ and from the proof of Lemma \ref{l:b} we deduce that the codimension of $\ul{\Delta}_i$ in $W_{\ul{K}(\ul{a}^{(i)})}(S)$ is at least
$i+2$. This implies that $\ul{\Delta}_i$ has codimension at least $2i+2$ in $S$. Since $S$ is of codimension $p$ in $\A^n$, we infer from the assumption $2i+2>n-p$ that $\ul{\Delta_i}$ is empty. Therefore $W_{\ul{K}(\ul{a}^{(i)})}(S)$ is smooth at
any point where $J\iFop$ has maximal rank.
\end{prf}

Proposition \ref{p:a} says that the lower dimensional generic polar varieties of $S$ (e.g. generic polar curves in case $p<n-1$, generic polar surfaces in case $p < n-2$ etc.) are empty or smooth at any point 
where $J\iFop$ has maximal rank.\spar

We are now going to show by an infinite family of examples that this conclusion is not always  true for higher dimensional generic polar varieties.\spar

\subsection{A family of singular generic polar varieties}\label{ss:3.1}

Let 
\[
n\ge 6,\;\; 
c:=\begin{bmatrix}
c_{1,1}&\cdots&c_{1,n}\\
c_{2,1}&\cdots&c_{2,n}
\end{bmatrix}\in \A^{2\times n}_{\R},\;\; 
a:=\begin{bmatrix}
a_{1,1}&\cdots&a_{1,n}\\
\vdots&\cdots&\vdots\\
a_{n-2,1}&\cdots&a_{n-2,n}\\
\end{bmatrix}\in \A^{(n-2)\times n}_{\R}
\]
such that the composed $(n\times n)$--matrix
$\scriptsize\begin{bmatrix}c\\a\end{bmatrix}$ represents a generic choice in $\A^{n\times n}$ which will be specified later. For the moment it suffices to suppose that 
$\scriptsize\begin{bmatrix}c\\a\end{bmatrix}$ and all $(2\times 2)$--submatrices of $c$ are regular and that all entries of $c$ are non--zero.
\spar

Let
\[
E_n:=\{(x_1\klk x_n)\in \A^n_{\R}\;|\;rk \begin{bmatrix}
c_{1,1}x_1&\cdots&c_{1,n}x_n\\
&a&
\end{bmatrix}=rk \begin{bmatrix}
c_{2,1}x_1&\cdots&c_{2,n}x_n\\
&a&
\end{bmatrix}=n-2\}
\]
and observe that $E_n$ is a linear subspace of $\A^n_{\R}$ of dimension at least $2(n-2)-n=n-4\ge 2$ (here $rk$ denotes the matrix rank).
\spar

Therefore we may assume without loss of generality that there exists an element $\xi=(\xi_1 \klk \xi_n)$ of $E_n$ with $\xi_1\neq 0$ and $\xi_2\neq 0$. Let
\[
c_1:=c_{1,1}\xi_1^2+\cdots +c_{1,n}\xi_n^2\;\;\text{and}\;\;c_2:=c_{2,1}\xi_1^2+\cdots +c_{2,n}\xi_n^2,
\]
\[
F^{(n)}_1:=c_{1,1}X_1^2+\cdots +c_{1,n}X_n^2-c_1,\;\;F^{(n)}_2:=c_{2,1}X_1^2+\cdots +c_{2,n}X_n^2-c_2
\]
and
\[
S^{(n)}:=\{F^{(n)}_1=F^{(n)}_2=0\}.
\]
For the sake of notational simplicity we shall write $E:=E_n$,  $F_1:=F^{(n)}_1$, 
$F_2:=F^{(n)}_2$ and $S:=S^{(n)}$.\spar

By construction $\xi$ belongs to $S_{\R}$, which implies that $S_{\R}$, and hence $S$, is non--empty.
\spar

Observe that any irreducible component of the complex variety $S$ has dimension at least $n-2\ge 4$.
\spar

Let $D$ be an arbitrary irreducible component of $S$. Thus we have $\dim\,D \ge 4$. We claim that $D$ contains a point $x=(x_1\klk x_n)$ such that there exists two indices $1\le u < v \le n$ with $x_u\neq 0$ and $x_v\neq 0$. Otherwise, for any $x=(x_1\klk x_n)$ of $D$ there exists an index $1\le j \le n$ with
$x_1=\cdots = x_{j-1}=x_{j+1}=\cdots =x_n=0$. From $0=F_1(x)=c_{1,j}x^2_j-c_1$
and $c_{1,j}\neq0$ we deduce that $S$, and therefore $D$, contains only finitely many such points. This implies $\dim\,D=0$ in contradiction to $\dim\,D\ge 4$. 
\spar

Therefore there exists a point $x=(x_1\klk x_n)$ of $D$ and indices $1\le u<v \le n$ with $x_u\neq0$ and $x_v\neq 0$. From 
\[J(F_1,\,F_2)= 2\,
\begin{bmatrix}
c_{1,1}X_1&\cdots&c_{1,n}X_n\\
c_{2,1}X_1&\cdots&c_{2,n}X_n\\
\end{bmatrix}
\]
and the ''genericity'' of $c$ we deduce that
\[
\det\begin{bmatrix}
\frac{\partial F_1}{\partial X_u}(x)&\frac{\partial F_1}{\partial X_v}(x)\\
\frac{\partial F_2}{\partial X_u}(x)&\frac{\partial F_2}{\partial X_v}(x)
\end{bmatrix}=4\,\det\begin{bmatrix}
c_{1,u}x_u&c_{1,v}x_v\\
c_{2,u}x_u&c_{2,v}x_v
\end{bmatrix}=4\,x_u x_v\,\det\begin{bmatrix}
c_{1,u}&c_{1,v}\\
c_{2,u}&c_{2,v}
\end{bmatrix}\neq 0
\]
holds. Thus the Jacobian $J(F_1,\,F_2)$ has maximal rank two
at the point $x$. Since $D$ was an arbitrary irreducible component of $S$, we conclude that any irreducible component of $S$ contains a point where $J(F_1,F_2)$ has maximal rank. This implies that the defining polynomials $F_1,\,F_2$ of $S$ form a reduced regular sequence in $\Q[X]$. In particular,
$S$ is of pure codimension two in $\A^n$. Let
\[
N_*(X):=\begin{bmatrix}
J(F_1,\,F_2)\\ a
\end{bmatrix}=\begin{bmatrix}
2c_{1,1}X_1&\cdots&2c_{1,n}X_n\\
2c_{2,1}X_1&\cdots&2c_{2,n}X_n\\
&a&
\end{bmatrix}.
\]
By construction $\xi$ belongs to $E$ and satisfies the conditions $F_1(\xi)=F_2(\xi)=0$. Therefore the real $((n-1)\times n)$--matrices
\[
\begin{bmatrix}
c_{1,1}\xi_1&\cdots&c_{1,n}\xi_n\\
&a&
\end{bmatrix}\;\;\text{and}\;\;
\begin{bmatrix}
c_{2,1}\xi_1&\cdots&c_{2,n}\xi_n\\
&a&
\end{bmatrix}\;\;\text{have rank}\;\;n-2.
\]
Moreover, the $((n-2)\times n)$--matrix $a$ has rank $n-2$. This implies that the vectors
$(c_{1,1}\xi_1\klk c_{1,n}\xi_n)$ and $(c_{2,1}\xi_1\klk c_{2,n}\xi_n)$ belong to the row span of $a$. Hence the real $(n\times n)$--matrix $N_{\ast}(\xi)$ has rank $n-2$.
Thus we have $\det N_{\ast}(\xi)=0$. 
\spar

Suppose for the moment that $a\in \A^{(n-2)\times n}$ and $\xi \in E$ were chosen in such a way that  $W_{\ul{K}^{n-3}(\ul{a})}(S)$ is, with respect to the properties of
polar varieties treated here, generic in the sense of Definition \ref{d:-1}.
From Observation \ref{o:a}
we deduce then that the polynomials $F_1,\,F_2$ and $\det N_{\ast}$ generate the ideal of definition of the polar variety $W_{\ul{K}^{n-3}(\ul{a})}(S)$. Hence $\xi$ belongs to the polar variety $W_{\ul{K}^{n-3}(\ul{a})}(S)$ which therefore turns out to be non--empty and of pure codimension three in $\A^n$. 
By the way, we see that $F_1,\,F_2$ and $\det N_{\ast}$ form a reduced 
regular sequence in $\Q[X]$.
\spar 

In the same vein as before we deduce from $\xi_1\neq 0$ and $\xi_2\neq 0$ that $\xi$ is a regular point of $S$. For $1\le u \le 2$ and $1\le l \le n$ 
we denote by $(-1)^{u+l}m_{u,l}$ the $(n-1)$--minor of 
$
\scriptsize\begin{bmatrix}
J(F_1,\,F_2)\\
a
\end{bmatrix}
$ obtained by deleting row number $u$ and column number $l$. From $\xi \in E$
we deduce $m_{u,l}(\xi)=0$.\spar
Let $1\le j \le n$. From the identity 
\[
\frac{\partial}{\partial X_j}\det N_{\ast}=
\sum_{1\le l \le n}( m_{1,l}\,\frac{\partial^2 F_2}{\partial X_l \,\partial X_j}+m_{2,l}\,\frac{\partial^2 F_1}{\partial X_l \,\partial X_j})=2(c_{2,j}m_{1,j}+c_{1,j}m_{2,j})
\]
we infer $(\frac{\partial}{\partial X_j}\det N_{\ast})(\xi)=0$. This implies
\[
\left((\frac{\partial}{\partial X_1}\det N_{\ast})(\xi)\klk (\frac{\partial}{\partial X_n}\det N_{\ast})(\xi)\right)=(0\klk 0).
\]
Therefore the rank of the Jacobian $J(F_1,\,F_2,\,\det N_{\ast})$ is two
at the point $\xi \in W_{\ul{K}^{n-3}(\ul{a})}(S)$.\spar

Since the polynomials $F_1,\,F_2$ and $\det N_{\ast}$ generate the ideal of definition of the polar variety $W_{\ul{K}^{n-3}(\ul{a})}(S)$, we conclude that 
$W_{\ul{K}^{n-3}(\ul{a})}(S)$ is singular at the point $\xi$.
\spar 

On the other hand, we deduce from $\xi_1\neq 0$ and $\xi_2\neq 0$ that
\[
\det\begin{bmatrix}
\frac{\partial F_1}{\partial X_1}(\xi)&\frac{\partial F_1}{\partial X_2}(\xi)\\
\frac{\partial F_2}{\partial X_1}(\xi)&\frac{\partial F_2}{\partial X_2}(\xi)
\end{bmatrix}
=4\,\xi_1\, \xi_2\,\det\begin{bmatrix}
c_{1,1}&c_{1,2}\\
c_{2,1}&c_{2,2}
\end{bmatrix}\neq 0
\]
holds. Thus the Jacobian $J(F_1, F_2)$ has maximal rank at $\xi$.
In particular, there exists a point where the polar variety $W_{\ul{K}^{n-3}(\ul{a})}(S)$ is not smooth and where $J(F_1, F_2)$ has maximal rank.\spar

We are now going to show that, for a suitable choice of $\xi$ in $E$, the polar variety  $W_{\ul{K}^{n-3}(\ul{a})}(S)$ is generic.\spar

For this purpose let 
\[
C=\begin{bmatrix}C_{1,1}&\cdots&C_{1,n}\\
C_{2,1}&\cdots&C_{2,n}\\
\end{bmatrix}
\]
and $(C_1,C_2)$ be a $(2\times n)$--matrix and a pair of new indeterminates. 
Further, let
 \[
A:=\begin{bmatrix}A_{1,1}&\cdots&A_{1,n}\\
&\cdots&\\
A_{n-2,1}&\cdots&A_{n-2,n}\\
\end{bmatrix}. 
 \] 
Observe that $S$ represents not a fixed algebraic variety but rather a constructible {\em family} 
of examples that depend on the choice of the $(n\times n)$--matrix $\scriptsize\begin{bmatrix}
c\\a\end{bmatrix}$ and the point $\xi \in E$ (recall that these data determine the values $c_1$ and $c_2$). To this family corresponds a genericity condition in the sense of 
Definition \ref{d:-1} which takes into account the properties of polar varieties we treat in this paper.  This genericity condition
may be expressed by the non-vanishing of a suitable non--zero polynomial $G\in \Q[A,C, C_1,C_2]$. The ''generic choice'' of the real $(n\times n)$--matrix 
$\scriptsize\begin{bmatrix}c\\a\end{bmatrix}$ means now that the bivariate polynomial $G(a,c,C_1,C_2)$ is non--zero and that all these matrices satisfy the requirements formulated at the beginning of this subsection.  
For the genericity of the polar variety $W_{\ul{K}^{n-3}(\ul{a})}(S)$ it suffices therefore to prove that the point $\xi=(\xi_1\klk \xi_n)$ may be chosen in $E$ in such a way that
\[
c_1:=c_{1,1}\xi_1^2+\cdots +c_{1,n}\xi_n^2\;\;\;\text{and}\;\;c_2:=c_{2,1}\xi_1^2+\cdots +c_{2,n}\xi_n^2
\]
satisfy the condition $G(a,c,c_1,c_2)\neq0$ (observe that $E$ depends only on the $(n\times n)$--matrix $\scriptsize\begin{bmatrix}c\\a\end{bmatrix}$ which we consider now as fixed).  To this end recall that $\dim\,E \ge 2$ holds and consider the polynomial map of affine spaces
$\varphi:E\to \A^2$ defined for $\xi=(\xi_1\klk \xi_n)$ in $E$ by 
\[\varphi(\xi):=(c_{1,1}\xi_1^2+\cdots +c_{1,n}\xi_n^2,\,c_{2,1}\xi_1^2+\cdots +c_{2,n}\xi_n^2).\]
It follows from our previous argumentation that for $\xi_1\neq0$ and $\xi_2\neq0$
the smooth map $\varphi$ is a submersion at $\xi$. Therefore the image of $\varphi$ is
Zariski dense in $\A^2$. This implies that there exists a point $\xi=(\xi_1\klk \xi_n)$ with
$\xi_1\neq0$ and $\xi_2\neq0$ in $E$ such that the bivariate polynomial $G(a,c,C_1,C_2)$
does not vanish at $(c_1,c_2):=\varphi(\xi)$. For such a choice of $\xi$ in $E$, 
the polar variety $W_{\ul{K}^{n-3}(\ul{a})}(S)$ turns out to be generic
in the sense stated before. \spar

In particular, we may suppose without loss of generality that
\[
c:=\begin{bmatrix}c_{1,1}&\cdots&c_{1,n}\\
c_{2,1}&\cdots&c_{2,n}\\
\end{bmatrix}\in \A^{2\times n}_{\R}\;\;\text{and}\;\;(c_1,c_2)\in \A^2_{\R}
\]
were chosen in such a way that all entries of $c$ and $c_1, c_2$ are non--zero, that all
$(2\times 2)$--submatrices of $c$ are regular and that, for any $1\le j \le n$, the values 
$\frac{c_1}{c_{1,j}}$ and $\frac{c_2}{c_{2,j}}$ are distinct.\spar 

Under this restriction the affine variety $S$ is empty or smooth.\spar

Assume that $S$ contains a singular point $x=(x_1\klk x_n)$. From our
preceding argumentation, we deduce that there exists an index $1\le k \le n$ such that
$x_j=0$ holds for any $1\le j\neq k \le n$. This implies $c_{1,k}x_k^2=c_1$
and $c_{2,k}x_k^2=c_2$ and therefore $\frac{c_1}{c_{1,k}}=\frac{c_2}{c_{2,k}}$
which contradicts our genericity condition. \spar

Therefore we have constructed for each $n\ge 6$ a smooth complete intersection variety
$S^{(n)}$ of pure codimension two and a generic polar variety $W_{\ul{K}^{n-3}(\ul{a})}(S^{(n)})$ of $S^{(n)}$ which contains a singular point.\spar

Observe that our construction may easily be modified in order to produce in the case
$p:=1$ {\em singular} generic {\em dual} polar varieties of regular hypersurfaces.
The behaviour of generic {\em classic} polar varieties is different in this case: they are always empty or smooth at regular points of the given hypersurface.\spar

Our family of  examples certifies that the {\em generic} higher dimensional polar  varieties of a smooth complete intersection variety $S$ may become singular. This contradicts the statement Theorem 10 (i) of \cite{bank4}.\spar

On the other hand we shall meet in Section \ref{s:4} three families of so--called 
{\em meagerly generic} polar varieties which are smooth and which therefore
satisfy the conclusion of \cite{bank4} Theorem 10 (i).\spar

Thus we face a quite puzzling situation, where the most general objects under consideration, namely the generic polar varieties, may become singular whereas
the more special meagerly generic polar varieties may turn out to be smooth.

\subsection{Desingularizing generic polar varieties}\label{ss:3.2}

For the rest of this section we suppose that $S$ is a smooth variety.
We are now going to try to explain in a more systematic way how singularities in higher
dimensional generic polar varieties of $S$ may arise.\spar

For expository reasons we limit our attention to the case of the classic polar varieties of $S$. 
Fix $1\le i \le n-p$ and consider the locally closed algebraic variety $\mathcal{H}_i$
of $\A^n\times \A^{(n-p-i+1)\times n}\times \P^{n-i}$ defined by
\begin{multline*}
\mathcal{H}_i:=\{(x,a,(\lambda:\vartheta))\in \A^n\times \A^{(n-p-i+1)\times n}\times \P^{n-i}\;|\;x\in S_{reg},\,rk\; a=n-p-i+1,\\
\lambda=(\lambda_1\klk \lambda_p)\in \A^p,\,\vartheta=(\vartheta_1\klk \vartheta_{n-p-i+1})\in \A^{n-p-i+1},\\ (\lambda,\vartheta)\neq 0,\;J\iFop(x)^T\cdot \lambda^T+a^T\cdot \vartheta^T=0\},
\end{multline*}
where $(\lambda:\vartheta):=(\lambda_1:\cdots :\lambda_p:\vartheta_1:\cdots : \vartheta_{n-p-i+1})$ belongs to $\P^{n-i}$, while  $(\lambda, \vartheta)\in \A^{n-i+1}$ is its affine counterpart and $J\iFop(x)^T$ denotes the transposed matrix of $J\iFop(x)$,
etc.\spar
Geometrically $\mathcal{H}_i$ may be interpreted as an incidence variety whose projection to the
space $\A^n \times \A^{(n-p-i+1)\times n}$ describes the locally closed variety
$\mathcal{W}_i$ of all pairs $(x,a)$ with $x\in S_{reg}$ and $a$ being a full--rank
matrix of $\A^{(n-p-i+1)\times n}$ such that $x$ belongs to $W_{\ul{K}(\ul{a})}(S)$. Below we shall see that $\mathcal{H}_i$ is non--empty, smooth and equidimensional. Thus $\mathcal{H}_i$ is a natural desingularisation of the variety
$\mathcal{W}_i$ in the sense of Room--Kempf (\cite{room, kem}).

\spar

Consider now an arbitrary point $(x,a,(\lambda:\vartheta))$ of $\mathcal{H}_i$ with
$a=[a_{k,l}]_{1\le k \le n-p-i+1 \atop{1\le l \le n}}$,
$\lambda=(\lambda_1\klk \lambda_p)\in \A^p$ and $\vartheta=(\vartheta_1\klk \vartheta_{n-p-i+1})\in \A^{n-p-i+1}$. Thus we have $(\lambda,\vartheta)\neq 0$.
We claim that $\vartheta \neq 0$ holds. Otherwise we would have $\vartheta =0$ and
$J\iFop(x)^T\cdot \lambda^T=0$. Since $x$ is a  regular point, this implies 
$\lambda=0$ in contradiction to $(\lambda,\vartheta)\neq 0$.\spar

Therefore we may assume without loss of generality $\vartheta_{n-p-i+1}=1$. Let $\tilde{a}:=[a_{k,l}]_{1\le k \le n-p-i \atop{1\le l \le n}}$ and denote for $1\le j \le n$  by $\theta_j$ 
the complex $(n\times (n-p-i+1))$--matrix containing
in row number $j$ the entries $1, \vartheta_1 \klk \vartheta_{n-p-i+1}$, whereas all other entries
of $\theta_j$ are zero. From $rk\, J\iFop(x)=p$ we deduce now that the complex $((n+p)\times (2n-i+(n-p-i+1)n))$--matrix 
\[
\begin{bmatrix}
J\iFop(x)&O_{p,p}&O_{p,n-p-i}&O_{p,n-p-i+1}&\cdots &O_{p,n-p-i+1}\\
\ast&J\iFop(x)^T&\tilde{a}^T& \theta_1&\cdots &\theta_n
\end{bmatrix}
\]
is of maximal rank $n+p$ (here $O_{p,p}$ denotes the $(p\times p)$--zero matrix, etc.). This implies that $\mathcal{H}_i$ is smooth and of dimension 
$n-p-i+(n-p-i+1)n$ at the point $(x,a,(\lambda:\vartheta))$. \\
In other words, the algebraic variety $\mathcal{H}_i$ is non--empty, smooth and equidimensional
of dimension $n-p-i+(n-p-i+1)n$. \spar

Let $\mu_i:\mathcal{H}_i\to \A^{(n-p-i+1)\times n}$ be the
canonical projection of $\mathcal{H}_i$ in $ \A^{(n-p-i+1)\times n}$ and suppose that for a generic choice of $ a\in \A^{(n-p-i+1)\times n} $ the polar variety $W_{\ul{K}^{n-p-i}(\ul{a})}(S)$ is non--empty.\spar

Observe that for any $x\in \A^n$ the following statements are equivalent
\begin{itemize}
\item[-] $x\in W_{\ul{K}^{n-p-i}(\ul{a})}(S)\cap S_{reg}$;
\item[-] there exist $\lambda\in \A^p$ and $\vartheta \in \A^{n-p-i+1}$ with $(\lambda,\vartheta)\neq0$ such that $x\in S_{reg}$ and $J\iFop(x)^T\cdot \lambda^T +a^T\cdot \vartheta^T=0$ holds;
\item[-]  there exist $\lambda\in \A^p$ and $\vartheta \in \A^{n-p-i+1}$ with $(\lambda,\vartheta)\neq0$ such that $(x,a,(\lambda:\vartheta))$ belongs to $\mu_i^{-1}(a)$.
\end{itemize}
Since by assumption $W_{\ul{K}^{n-p-i}(\ul{a})}(S)$ is non--empty for any generic choice $a\in \A^{(n-p-i+1)\times n}$,
we conclude that
$\mu_i^{-1}(a)\neq \emptyset$ holds. 
This implies that the image of the morphism of equidimensional algebraic varieties $\mu_i:\mathcal{H}_i\to \A^{(n-p-i+1)\times n}$ is Zariski dense in $\A^{(n-p-i+1)\times n}$.\spar

Let $a=[a_{k,l}]_{1\le k \le n-p-i+1 \atop{1\le l \le n}}$ be a fixed generic choice in $ \A^{(n-p-i+1)\times n}$ and let 
\begin{multline*}
H_i:=\{(x,(\lambda:\vartheta))\in \A^n\times  \P^{n-i}\;|\;x\in S_{reg},\,
\lambda\in \A^p,\,\vartheta \in \A^{n-p-i+1},\\ (\lambda,\vartheta)\neq 0,\;J\iFop(x)^T \cdot\lambda^T+a^T\cdot \vartheta^T=0\},
\end{multline*}
Since the morphism $\mu_i$ has in $\A^{(n-p-i+1)\times n}$ a Zariski dense image, $a$ is chosen
generically in $\A^{(n-p-i+1)\times n}$ and $\mathcal{H}_i$ is 
smooth and equidimensional of dimension $n-p-i+(n-p-i+1)n$, we deduce from the Weak Transversality Theorem of Thom--Sard that $H_i$ is
non--empty, smooth and equidimensional of dimension $n-p-i$.\\
Let $\omega_i:H_i\to \A^n$ be the canonical projection $H_i$ in $\A^n$. From our previous reasoning we conclude that the image of $\omega_i$ is $W_{\ul{K}^{n-p-i}(\ul{a})}(S)\cap S_{reg}$. Since this algebraic variety is equidimensional and of the same dimension as $H_i$, namely $n-p-i$,
we infer that there exists a Zariski dense open subset $U$ of $W_{\ul{K}^{n-p-i}(\ul{a})}(S)$ such that the $\omega_i$--fiber of any $x\in U$ is zero--dimensional.\spar

Using the notations of Section \ref{s:3} let us now consider an arbitrary point\\ $x\in W_{\ul{K}^{n-p-i}(\ul{a})}(S)\cap S_{reg}$ which does not belong to $\ul{\Delta}_i$. Without loss of generality we may assume that
\[
m:=\det\;\begin{bmatrix}
\frac{\partial F_1}{\partial X_1}&\cdots&\frac{\partial F_1}{\partial X_{n-i}}\\
\vdots&\cdots&\vdots\\
\frac{\partial F_p}{\partial X_1}&\cdots&\frac{\partial F_p}{\partial X_{n-i}}\\
a_{1,1}&\cdots&a_{1,n-i}\\
\vdots&\cdots&\vdots\\
a_{n-p-i,1}&\cdots&a_{n-p-i,n-i}\\
\end{bmatrix}.
\]
does not vanish at $x$.\\
Let $\lambda\in A^p$ and $\vartheta \in \A^{n-p-i+1}$ with $\vartheta=(\vartheta_1\klk \vartheta_{n-p-i+1})$ and $(\lambda,\vartheta)\neq 0$ such that $(x,(\lambda:\vartheta))$ belongs to $\omega_i^{-1}(x)$. We claim $\vartheta_{n-p-i+1}\neq0$. Otherwise we would have for $\tilde{a}:=[a_{k,l}]_{1\le k \le n-p-i \atop{1\le l \le n}}$ and $\tilde{\vartheta}=(\vartheta_1\klk \vartheta_{n-p-i})$ that $(\lambda, \tilde{\vartheta})\neq0$ and $J\iFop(x)^T \cdot \lambda^T+
\tilde{a}^T\cdot \tilde{\vartheta}^T=0$ holds. This in turn would imply that the rank of the complex $((n-i)\times n)$--matrix $\scriptsize\begin{bmatrix}
J\iFop(x)\\ \tilde{a}
\end{bmatrix}$ is at most $n-i-1$, whence $m(x)=0$, a contradiction.\spar

Therefore we may assume without loss of generality $\vartheta_{n-p-i+1}=1$. Since the algebraic variety $H_i$ is smooth at the point $(x,(\lambda:\vartheta))$, we conclude that the complex $((n+p)\times (2n-i))$--matrix 
\[
\begin{bmatrix}
J\iFop(x)&O_{p,p}&O_{p,n-p-i}\\
\ast&J\iFop(x)^T&\tilde{a}^T
\end{bmatrix}
\]
is of maximal rank $n+p$. From $m(x)\neq0$ we deduce that the submatrix\\  $\scriptsize\begin{bmatrix}J\iFop(x)^T&\tilde{a}^T \end{bmatrix}$ is of maximal rank $n-i$. Taking into account
that\\$x\in (W_{\ul{K}^{n-p-i}(\ul{a})}(S)\cap S_{reg})\setminus \ul{\Delta}_i$
is a smooth point of $W_{\ul{K}^{n-p-i}(\ul{a})}(S)\cap S_{reg}$, we infer that there
exists in the Euclidean topology a sufficiently small open neighborhood $O$ of $x$ in $W_{\ul{K}^{n-p-i}(\ul{a})}(S)\cap S_{reg}$ and a smooth map $\mu:O\to H_i$ with $\mu(x)=
(x,(\lambda:\vartheta))$ and $\omega_i \circ \mu= id_{\,O}$, where $ id_{\,O}$
denotes the identity map on $O$. In particular, 
${\omega_i}_{|_{\omega_i^{-1}(O)}}\!:\omega_i^{-1}(O) \to U$ 
is a submersion at $(x,(\lambda:\vartheta))$.
Translating this reasoning to the language of commutative algebra, we conclude that $\omega_i$
as morphism of algebraic varieties $H_i\to W_{\ul{K}^{n-p-i}(\ul{a})}(S)$ is unramified at any point of $\omega_i^{-1}(x)$. In particular, $\omega_i^{-1}(x)$ is a zero--dimensional algebraic variety.\spar

Let us now consider an arbitrary point $x\in \ul{\Delta}_i\cap S_{reg}$.  Then the complex $((n-i+1)\times n)$--matrix $\scriptsize\begin{bmatrix}J\iFop(x)\\a\end{bmatrix}$ has rank at most $n-i-1$. Hence the linear variety
\begin{multline*}
\Gamma:=\{(\lambda:\vartheta)\in \P^{n-i}\;|\;\lambda \in \A^p,\, \vartheta \in \A^{n-p-i+1},\,
(\lambda,\vartheta)\neq0, \\J\iFop(x)^T\cdot \lambda^T+a^T\cdot \vartheta^T=0\},
\end{multline*}
which is isomorphic to $\omega_i^{-1}(x)$, has dimension at least one. Thus $\omega_i^{-1}(x)$ is not a zero--dimensional algebraic variety. In particular, if  
$W_{\ul{K}^{n-p-i}(\ul{a})}(S)$ is smooth at $x$, then for any sufficiently  small open neighborhood $O$ of $x$ in
$W_{\ul{K}^{n-p-i}(\ul{a})}(S)$ and any $\lambda \in \A^p,\;\vartheta \in \A^{n-p-i+1}$ with $(\lambda, \vartheta)\neq0$ and $(x,(\lambda : \vartheta))\in \omega_i^{-1}(x)$ the smooth map ${\omega_i}_{|_{\omega_i^{-1}}(O)}:\omega_i^{-1}(O) \to O$ is not a submersion at
$(x,(\lambda:\vartheta))$.\spar

We may summarize this argumentation by the following statement.

\begin{proposition}\label{p:b}
Let notations be as above and suppose that for $a\in \A^{(n-p-i+1)\times n}$ generic the polar variety
$W_{\ul{K}^{n-p-i}(\ul{a})}(S)$ is not empty. Then we have for any regular
point $x$ of $W_{\ul{K}^{n-p-i}(\ul{a})}(S)$:
$x$ belongs to $\ul{\Delta}_i$ if and only if $\omega_i^{-1}(x)$ is a zero-dimensional variety. If this is the case, the morphism of algebraic varieties $\omega_i:H_i\to W_{\ul{K}^{n-p-i}(\ul{a})}(S) $ is unramified at any point  of  $\omega_i^{-1}(x)$.
\end{proposition}

The morphism of algebraic varieties $\omega_i:H_i\to W_{\ul{K}^{n-p-i}(\ul{a})}(S) $ represents a natural desingularization of the open generic polar variety 
$W_{\ul{K}^{n-p-i}(\ul{a})}(S)\cap S_{reg}$ in the spirit of Room--Kempf  \cite{room, kem}. Proposition \ref{p:b} characterizes $\ul{\Delta}_i$ as a kind of an algebraic--geometric
''discriminant locus'' of the morphism $\omega_i$ which maps the smooth variety $H_i$ to the possibly non--smooth generic polar variety $W_{\ul{K}^{n-p-i}(\ul{a})}(S)$.\spar

We finish this section by a decomposition of the locus of the regular points
of the generic classic polar varieties of $S$ into smooth Thom--Boardman strata  
of suitable polynomial maps defined on $S_{reg}$.
Our main tool will be Mather's theorem on generic projections \cite{math} (see
also \cite{ottav}). Then we will discuss this decomposition in the light
of Proposition \ref{p:a}.\spar

Let us fix an index $1 \le i \le n-p$ and a complex $(n\times (n-p-i+1))$--matrix 
$a\in \A^{(n-p-i+1) \times n}$ of rank $n-p-i+1$. Then $a$ induces a linear
map $\Pi_a:\A^n\to \A^{n-p-i+1}$ which is for $x\in \A^n$ defined by $x\,a$.
Let $\pi_a$ be the restriction of this linear map to the variety $S$.
\spar

Let $x$ be an arbitrary point of $S$. We denote by $T_x:=T_x(S)$ the Zariski tangent space of $S$ at $x$, by $d_x(\pi_a): T_x \to \A^{n-p-i+1}$ the corresponding linear map
of tangent spaces at $x$ and $\pi_a(x)$, by $ker\,d_x(\pi_a)$ its kernel and by $\dim_{\C}(ker\,d_x(\pi_a))$ the dimension of the $\C$--vector space $ker\,d_x(\pi_a)$. \spar

Furthermore,
for $i\le j \le n-p$ we denote by
\[
\Sigma^{(j)}(\pi_a):=\{x\in S_{reg}\;|\;\dim_{\C}(ker\,d_x(\pi_a))=j\}
\]
the corresponding Thom--Boardman stratum of $\pi_a$. Observe that $\Sigma^{(j)}(\pi_a)$ is a constructible subset of $\A^n$. With these notations we have the following result.
\begin{lemma}\label{l:d}
\begin{itemize}
\item [(i)]$W_{\ul{K}(\ul{a})}(S)\cap S_{reg}=\bigcup_{i\le j\le n-p}\Sigma^{(j)}(\pi_a)$.
\item [(ii)] $\ul{\Delta_i} \cap S_{reg}=\bigcup_{i+1\le j\le n-p}\Sigma^{(j)}(\pi_a)$.
\end{itemize}
\end{lemma}
\begin{prf}
Let $x$ be an arbitrary point of $S_{reg}$ and let $N_x$ and $E_a$ be the $\C$--vector spaces generated by the rows of the matrices $J\iFop(x)$ and $a$, respectively. Thus we have $\dim_{\C}\,T_x=n-p,\; \dim_{\C}\,N_x=p$ and $\dim_{\C}\,E_a=n-p-i+1$.\spar
Let $\yon$ be new indeterminates. With respect to the bilinear form $X_1Y_1+\cdots +X_nY_n$ we may consider $T_x$ and $N_x$ as orthogonal subspaces of $\A^n$.\spar

From the definition of the polar variety $W_{\ul{K}(\ul{a})}(S)$ we deduce 
that the point $x$ belongs to it if and only if $\dim_{\C}\,(N_x+E_a)\le n-i$ holds. 
Taking into account the identity
\[
\dim_{\C}\,(N_x+E_a)= \dim_{\C}\,N_x+\dim_{\C}\,E_a- \dim_{\C}\,(N_x\cap E_a)
\]
we conclude that the condition $\dim_{\C}\,(N_x+E_a)\le n-i$ is equivalent to
\[
\dim_{\C}\,(N_x\cap E_a)\ge \dim_{\C}\,N_x+\dim_{\C}\,E_a - n+i =p+(n-p-i+1)-n+i=1.
\]
Since the $\C$--linear spaces $N_x\cap E_a$ and $T_x$ are orthogonal subspaces of $\A^n$ with 
respect to the bilinear form $X_1Y_1+\cdots +X_nY_n$, we infer 
\[
\dim_{\C}\,\Pi_a(T_x)\le \dim_{\C}\,E_a - \dim_{\C}\,(N_x\cap E_a)=n-p-i+1 - \dim_{\C}\,(N_x\cap E_a).
\]
Thus the conditions
\[
\dim_{\C}(ker\,d_x(\pi_a))\ge i,\;\;\dim_{\C}\,\Pi_a(T_x)\le n-p-i\;\;\text{and}\;\;
x\in  W_{\ul{K}(\ul{a})}(S)\cap S_{reg}
\]
are equivalent. This proves the statement $(i)$ of the lemma.\spar

Statement $(ii)$ can be shown similarly, observing that the condition $x\in \ul{\Delta_i} \cap S_{reg}$ is equivalent to $\dim_{\C}\,(N_x + E_a)\le n-i-1$.
\end{prf}

The following deep result represents the main ingredient of our analysis.\spar
\begin{lemma}\label{l:e}
For an integer $j\ge i$, let $v_{i}(j):=(j-i+1)j$. There exists a non--empty Zariski open 
subset $U_i$ of complex full--rank $((n-p-i+1)\times n)$--matrices of $\A^{(n-p-i+1)\times n}$ 
such that for any any $a\in U_i$ and any index $i\le j \le n-p$ the constructible set
$\Sigma^{(j)}(\pi_a)$ is empty or a smooth subvariety of $S$ of codimension $v_{i}(j)$.
\end{lemma}
The proof of Lemma \ref{l:e} is based on two technical statements which represent particular instances of  \cite{ottav}, Theorem 2.31 on one side and Theorems 3.10
and 3.11 on the other. We will not reproduce the details here. These statements go back to Mather's original work \cite{math}.\spar

Combining Lemma \ref{l:d} and Lemma \ref{l:e} we obtain for generic $a\in \A^{(n-p-i)\times n}$ decompositions of  $W_{\ul{K}(\ul{a})}(S)\cap S_{reg}$ and
$\ul{\Delta_i} \cap S_{reg}$ in pieces which are empty or smooth and of predetermined codimension in $S$.\spar

In order to illustrate the power of this analysis we are now giving an alternative proof of
Proposition \ref{p:a} using the following argumentation.\spar

Suppose $2i+2> n-p$ and let $a$ be a complex full--rank $((n-p-i+1)\times n)$--matrix
contained in the set $U_i$ introduced in Lemma \ref{l:e}.\spar

Then we have for $i< j \le n-p$
\[
n-p-v_{i}(j)\le n-p-v_{i}(i+1)=n-p-(2i+2)<0
\]
and therefore
\[\Sigma^{(j)}(\pi_a)=\emptyset.\]
Lemma \ref{l:d} implies now 
\[W_{\ul{K}(\ul{a})}(S)\cap S_{reg}=\Sigma^{(i)}(\pi_a)\]
and from Lemma \ref{l:e} we deduce finally that $W_{\ul{K}(\ul{a})}(S)\cap S_{reg}$ is empty or a smooth subvariety of $S$ of codimension $v_{i}(i)=i$.\spar

\section{Meagerly generic polar varieties}\label{s:4}

In terms of the notions and notations introduced in Section \ref{s:1} and Section \ref{s:2}, for $a:=\left[ a_{k,l}\right]_{{1 \le k \le n-p} \atop {1 \le l \le n}}$ 
being a complex full
rank matrix, the (geometric) degrees of the classic and dual polar varieties
$W_{\ul{K}(\ul{a}^{(i)})}(S)$ and $W_{\ol{K}(\ol{a}^{(i)})}(S)$, $1 \le i \le n-p$,
are constructible functions of $a\in \A^{(n-p)\times n}$. In view of Definition \ref{d:-1}, this implies that, for $a\in \A^{(n-p)\times n}$ generic, these degrees are independent of $a$ and therefore invariants of $S$. We denote by ${\delta}_{classic}$
the maximal geometric degree of the classic polar varieties $W_{\ul{K}(\ul{a}^{(i)})}(S)$, $1 \le i \le n-p$, when $a\in \A^{(n-p)\times n}$ is generic. 
In a similar way we define ${\delta}_{dual}$ as the maximal degree of the
$W_{\ol{K}(\ol{a}^{(i)})}(S)$, $1 \le i \le n-p$, for $a\in \A^{(n-p)\times n}$ generic.
Thus ${\delta}_{classic}$ and ${\delta}_{dual}$ are also well--defined invariants of $S$. The main outcome of Theorem \ref{th:c} below is the observation that   ${\delta}_{classic}$
and ${\delta}_{dual}$ dominate (in a suitable sense) the degrees of {\em all}
classic and dual polar varieties of $S$.\spar

The statement of Theorem \ref{th:c} is motivated by our aim to apply the geometry
of polar varieties to algorithmic problem of real root finding. This can be explained as follows:\spar
Suppose that there is given a division--free arithmetic circuit $\beta$ in
$\Q[X]$ evaluating the reduced regular sequence $\Fop \in \Q[X]$ that
defines the complex variety $S$. Let $L$ be the size of $\beta$ (which
measures the number of arithmetic operations required by $\beta$ for the
evaluation of $\Fop$) and let $d$ and $D$ be the maximal degrees of the
polynomials $\Fop$ and the intermediate varieties $\{F_1=0\klk F_k=0\}, 1\le k\le p$,
respectively.
Suppose furthermore that the real trace $S_{\R}$ of
$S$ is non--empty, smooth and of pure codimension $p$ in $\A_{\R}^n$.
\spar

Then a set of at most ${\delta}_{dual}$ real algebraic sample points for
the connected components of $S_{\R}$ can be found using at most 
$O(\binom{n}{p}L n^4 p^2d^2(\max\{D,\delta_{dual}\})^2)$ arithmetic operations and
comparisons in $\Q$ (see \cite{bank3, bank4}).\\

If additionally $S_{\R}$ is compact, the same conclusion holds true with
${\delta}_{classic}$ instead of ${\delta}_{dual}$.\\
The underlying algorithm can be implemented in the {\em non--uniform deterministic} or
alternatively in the {\em uniform probabilistic} complexity model.
\spar

From the B\'ezout Inequality \cite{he} one deduces easily that
${\delta}_{classic}$ and ${\delta}_{dual}$ are bounded by ${d^n}{p^{n-p}}$
(see \cite {bank3} and \cite {bank4} for proofs).
This implies that from the geometrically unstructured and extrinsic point
of view, the quantities ${\delta}_{classic}$ and ${\delta}_{dual}$ and hence the
complexity of the procedures of \cite {bank1}, \cite {bank2}, \cite
{bank3}, \cite {bank4}, \cite {mohab2} and \cite{SaSch} are in worst case of order
$(nd)\,^{O(n)}$. This meets all previously known algorithmic bounds (see e.g.
\cite{grivo}, \cite{HeRoSo}, \cite {baporo}, \cite{RoRoSa}, \cite{Sa05}). 
From \cite{GiHe01} and \cite{CaGiHeMaPa} one deduces a {\em worst case} lower
bound of order $d^{\Omega(n)}$ for the complexity of the real root finding problem
under consideration. Thus complexity improvements are only possible if we distinguish
between well-- and ill--posed input systems $F_1=0\klk F_p=0$ or varieties $S$. We realize this distinction by means of the invariants ${\delta}_{classic}$ and ${\delta}_{dual}$.
\spar

For details on the notion of the (geometric) degree of algebraic varieties we
refer to \cite{he} (or \cite{fu} and \cite{vo}). For expository
reasons we define the degree of the empty set as zero.
For algorithmic aspects we refer to \cite{bank2}, \cite{bank3}, \cite{bank4}, \cite {mohab2} and \cite {SaSch}.
The special case of $p:=1$ and $S_{\R}$ compact is treated in \cite{bank1}.
\spar

We mention here that the algorithms in the papers \cite{bank1} and \cite{bank2},
which treat the special case $p:=1$ and the general case $1 \le p \le n$ 
when $S_{\R}$ is compact, use instead of the {\em generic}
polar varieties $W_{\ul{K}(\ul{a}^{(i)})}(S)$, $1 \le i \le n-p$ ,
with $a:=\left[ a_{k,l}\right]_{{1 \le k \le n-p} \atop {1 \le l \le n}}$
generic, more special polar varieties, where 
$a:=\left[ a_{k,l}\right]_{{1 \le k \le n-p} \atop {1 \le l \le n}}$ ranges over a (Zariski) dense set of rational points of a suitable irreducible subvariety of $\R^{(n-p) \times n}$.
In the sequel we shall call these special polar varieties {\em meagerly
generic} (see Definition \ref{d:1} below).
\spar

It is not difficult to show that the argumentation of \cite{bank1} in the case
$p=1$ covers also the instance of generic polar varieties. 
However, in the case $1< p \le n$, the polar varieties
$W_{\ul{K}(\ul{a}^{(i)})}(S)$, $1 < i \le n-p$ appearing in
\cite{bank2} are not generic, though they turn out to be
non--empty, smooth and of pure codimension $i$ in $S$ (see Example 1).
\spar

Although it was not essential for the final outcome of \cite{bank2}, namely the 
correctness and complexity bound of the proposed point finding  algorithm, we
believed at the moment of publishing this paper that the smoothness of
these meagerly generic polar varietes implies the smoothness of their generic counterparts. As we have seen in Section \ref{s:3}, this
conclusion is wrong. 
\spar

In order to estimate the complexity of root finding algorithms based on meagerly generic polar varieties (as in \cite{bank2}) in terms of not too artificial invariants of the underlying variety $S$ (like ${\delta}_{classic}$ or ${\delta}_{dual}$), we need to know that the degrees of the generic polar varieties dominate the degrees of their more special meagerly generic counterparts. This is the aim of the next result.\spar 

\begin{theorem}\label{th:c}
Let $a:=\left[ a_{k,l}\right]_{{1 \le k \le n-p} \atop {1 \le l \le n}}$
 be a complex $((n-p) \times n)$--matrix. Suppose that for $1 \le i \le n-p$
the $((n-p-i+1) \times n)$--matrix $\left[ a_{k,l}\right]_{{1 \le k \le
n-p-i+1} \atop {1 \le l \le n}}$ is of maximal rank $n-p-i+1$ and that the
polar varieties $W_{\ul{K}(\ul{a}^{(i)})}(S)$ and
$W_{\ol{K}(\ol{a}^{(i)})}(S)$ are empty or of pure codimension $i$
in $S$.\spar
Furthermore, let be given a non--empty Zariski open set $U$ of complex
$((n-p) \times n)$--matrices $b:=\left[ b_{k,l}\right]_{{1 \le k \le n-p}
\atop {1 \le l \le n}}$ such that for each $b \in U$ and each $1 \le i \le
n-p$ the $((n-p-i+1) \times n)$--matrix  $\left[ b_{k,l}\right]_{{1 \le k
\le n-p-i+1} \atop {1 \le l \le n}}$ is of maximal rank $n-p-i+1$ and such
that the polar varieties $W_{\ul{K}(\ul{b}^{(i)})}(S)$ and
$W_{\ol{K}(\ol{b}^{(i)})}(S)$ are uniformly empty or of
codimension $i$ in $S$ and have fixed geometric degrees, say
${{\delta}_{classic}}^{(i)}$ and ${{\delta}_{dual}}^{(i)}$.\spar

Suppose finally that, for $1 \le i \le n-p$, the non--emptiness of
$W_{\ul{K}(\ul{a}^{(i)})}(S)$ (or
$W_{\ol{K}(\ol{a}^{(i)})}(S)$) implies that of
$W_{\ul{K}(\ul{b}^{(i)})}(S)$ (or
$W_{\ol{K}(\ol{b}^{(i)})}(S)$).
\spar
Then we have for $1 \le i \le n-p$ the estimates
\[\deg\,W_{\ul{K}(\ul{a}^{(i)})}(S)
 \le {{\delta}_{classic}}^{(i)}\]
and
\[\deg\,W_{\ol{K}(\ol{a}^{(i)})}(S)
 \le {{\delta}_{dual}}^{(i)}.\]
\end{theorem}

The proof of Theorem \ref{th:c} is based on Lemma \ref{l:c} below. 
In order to state this lemma, suppose that there is given a morphism of affine varieties $\varphi: V\to \A^r$, where $V$ is equidimensional of dimension $r$
and the restriction map ${\varphi}_{|_C}: C\to \A^r$ is dominating for each irreducible component $C$ of $V.$ Then $\varphi$ induces an embedding of the rational function field $\C(\A^r)$ of $\A^r$ into the total quotient ring $\C(V)$
of $V$. Observe that $\C[V]$ is a subring of $\C(V)$ and that $\C(V)$ is isomorphic to 
the direct product of the function fields of all irreducible components of $V$. Hence 
$\C(V)$ is a finite dimensional $\C(\A^r)$--vector space and
the degree of the morphism $\varphi$, namely 
$\deg\,\varphi:= \dim_{\C(\A^r)}\,\C(V)$, is a well--defined quantity. 
Finally, we denote for $y\in \A^r$
by $\#\, \varphi^{-1}(y)$ the cardinality of the $\varphi$--fiber $\varphi^{-1}(y)$ at
$y$ and by $\#_{isolated}\, \varphi^{-1}(y)$ the (finite) number of its isolated points. 
\spar

With these assumptions and notations we are now going to show the following
elementary geometric fact which represents a straightforward generalization
of \cite{he} Proposition 1 and its proof.
\begin{lemma}\label{l:c}
\begin{itemize}
\item [(i)] Any point $y\in \A^r$ satisfies the condition
\[
\#_{isolated}\, \varphi^{-1}(y)\le \deg\, \varphi.\]
In particular, if the fiber of $\varphi$ is finite, we have $\#\, \varphi^{-1}(y)\le \deg\, \varphi$.
\item [(ii)] There exists a non--empty Zariski open subset $U\subset \A^r$ such that
any point $y\in U$ satisfies the condition
\[
\#\, \varphi^{-1}(y) = \deg\, \varphi.
\]
\end{itemize}
\end{lemma}
\begin{prf}
Let us show the first statement of Lemma \ref{l:c}$(i)$, the second one is then obvious. We proceed by induction on $r$.\spar
The case $r=0$ is evident. Let us therefore suppose $r>0$ and that the lemma is true
for any morphism of affine varieties $\varphi': V'\to \A^{r-1}$ which satisfies the previous requirements.\spar
Consider an arbitrary point $y=(y_1\klk y_r)$ of $\A^r$.
If $\varphi^{-1}(y)$ does not contain isolated points we are done. Therefore we may
consider an arbitrary isolated point $z$ of $\varphi^{-1}(y)$. Denoting the coordinate functions of
$\A^r$ by $Y_1 \klk Y_r$, let us identify the hyperplane $\{Y_1-y_1=0\}$ of $\A^r$ with $\A^{r-1}$. Thus we have $y\in \A^{r-1}$. From $\varphi^{-1}(y)\neq \emptyset$ and
the Dimension Theorem for algebraic varieties (see e.g. \cite{sha}, Vol. I, Ch. I, Corollary
1 of Theorem 6.2.5) we deduce that $\varphi^{-1}(\A^{r-1})$ is a non--empty subvariety of
$V$ of pure codimension one. Since $y$ belongs to $\A^{r-1}$ and $\varphi^{-1}(y)$
contains an isolated point we deduce from the Theorem on Fibers (see e.g. \cite{sha}, Vol. I, Ch. I, Theorem 6.3.7) that there exists an irreducible component $C'$ of
$\varphi^{-1}(\A^{r-1})$ such that the induced morphism of irreducible affine varieties 
$\varphi_{|_{\C'}}: C'\to \A^{r-1}$ is dominating.\spar
Let $C''$ be an irreducible component of $\varphi^{-1}(\A^{r-1})$ such that $\varphi(C'')$  is
not Zariski dense in $\A^{r-1}$ and denote the Zariski closure of $\varphi(C'')$ in $\A^{r-1}$ by $D$. Then $D$ is an irreducible subvariety of $\A^{r-1}$ of dimension at most $r-2$.\spar
We claim now that $z$ does not belong to $C''$. Otherwise, by the Theorem on Fibers, all components of the 
$\varphi_{|_{\C''}}$--fiber of $y$, namely $\varphi^{-1}(y)\cap C''$,
have dimension at least $\dim\,C'' - \dim\,D > 0$ and
$z$ belongs to $\varphi^{-1}(y)\cap C''$. This contradicts our assumption that $z$
is an isolated point of  $\varphi^{-1}(y)$.\spar

Therefore the isolated points of $\varphi^{-1}(y)$ are contained in $V':= \bigcup_{C'\in \mathcal{Z}} \;C'$, where
\[ \mathcal{Z}:=\{C'\,|\,C'\;\text{irreducible component of}\;\varphi^{-1}(\A^{r-1}),\;\varphi_{|_{\C'}}: C'\to \A^{r-1}\;\text{dominating}\}.\]
Observe that $V'$ is a (non--empty) equidimensional affine variety of dimension $r-1$ and that $\varphi$
induces a morphism of affine varieties $\varphi':V'\to \A^{r-1}$ such that $\varphi'_{|_{C'}}:C'\to \A^{r-1}$ is dominating for any irreducible component $C'$ of $V'$. Moreover, $y$ belongs to the image of $\varphi'$ and 
$(\varphi')^{-1}(y)$ and $\varphi^{-1}(y)$ have the same isolated points.\spar

Hence, by induction hypothesis, the statement $(i)$ of the lemma is valid for $\varphi'$ and   
$y$. Consequently we have 
\[
\#_{isolated}\,\varphi^{-1}(y)=\#_{isolated}\,(\varphi')^{-1}(y)\le \deg\,\varphi'.
\]
We finish now the proof of lemma \ref{l:c},$(i)$ by showing that $\deg\,\varphi' \le \deg\,\varphi$
holds.\spar

For this purpose let us abbreviate $Y:=(Y_1\klk Y_r)$ and let us write $\C[Y]$ and
$\C(Y)$ instead of $\C[\A^r]$ and $\C(\A^r),$ respectively. Thus we consider $\C[Y]$
as a $\C$--subalgebra of $\C[V]$ and $\C(Y)$ as a $\C$--subfield of the total quotient ring
$\C(V)$.\spar

For any irreducible component $C$ of $V$ the morphism of affine varieties $\varphi_{|_{C}}:C \to \A^r$ is by assumption dominating. This implies that $Y_1 - y_1$ is not
a zero divisor of $\C[V]$ or $\C(V)$.\spar

The canonical embedding of $V'$ in $V$ induces a surjective $\C$--algebra homomorphism $\C[V]\to \C[V']$ which associates to each $g\in \C[V]$	an image
in $\C[V']$ denoted by $g'$. Therefore there exists for $\delta:=\deg\,\varphi'$ 	
elements $g_1\klk g_{\delta}$ of $\C[V]$ such that $g'_1\klk  g'_{\delta}$
form a vector space basis of $\C(V')$ over $\C(\A^{r-1})=\C(Y_2 \klk Y_r)$.
In particular, $g'_1\klk  g'_{\delta}$ are $\C(Y_2 \klk Y_r)$--linearly independent.
We claim that $g_1\klk g_{\delta}$, considered as elements of $\C(V)$, are also linearly
independent over $\C(\A^r)=\C(Y).$\\
Otherwise there exist polynomials $b_1\klk b_{\delta}\in \C[Y]$ not all zero, satisfying the condition
\begin{equation}\label{eq:3}
b_1\,g_1+\cdots + b_{\delta}\,g_{\delta}=0
\end{equation}
in $\C[V]$.\spar
For $1\le j \le \delta$ let $b_j=(Y_1-y_1)^{u_j}c_j$, where $u_j$ is a non--negative integer and $c_j$ is polynomial of $\C[Y]$ which is not divisible by $Y_1-y_1$. Thus (\ref{eq:3}) can be written as
\[
(Y_1-y_1)^{u_1}c_1\,g_1+ \cdots + (Y_1-y_1)^{u_{\delta}}c_{\delta}\,g_{\delta}=0.
\]
Therefore, since $Y_1-y_1$ is not a zero divisor in $\C[V]$, we may assume without loss of generality that there exists $1\le j_0 \le \delta$ with $u_{j_0}=0$. This means that we may suppose $b_{j_0}$ is not divisible by $Y_1-y_1$ in $\C[Y]$.\spar

Condition (\ref{eq:3}) implies now that 
\begin{equation}\label{eq:4}
b'_1\,g'_1+ \cdots + b'_{\delta}\,g'_{\delta}=0
\end{equation}
holds in $\C[V']$ and therefore also in $\C(V')$. Observe that $b'_j$ is the polynomial of $\C[Y_2\klk Y_r]$ obtained by substituting in $b_j$ the indeterminate $Y_1$ by $y_1\in \C$, namely $b'_j=b_j(y_1, Y_2\klk Y_r)$, where $1\le j \le \delta$.   
Since by assumption $b_{j_0}$ is not divisible by $Y_1-y_1$ in $\C[Y]$, we have 
$b'_{j_0}=b_{j_0}(y_1, Y_2\klk Y_r)\neq 0$. Thus (\ref{eq:4}) implies that 
$g'_1\klk  g'_{\delta}\in \C(V')$ are linearly dependent over $\C(Y_2\klk Y_r)=\C(\A^{r-1})$ which contradicts the hypothesis that $g'_1\klk  g'_{\delta}$ is a vector space basis of $\C(V')$.\spar

Thus $g_1\klk  g_{\delta}\in \C(V)$ are linearly independent over $\C(\A^r)$. Hence we may conclude 
\[
\deg\,\varphi'=\dim_{\,\C(\A^{r-1})}\,\C(V')=\delta \le  \dim_{\,\C(\A^r)}\,\C(V)=\deg\,\varphi
\]
and we are done.\spar

We are now going to prove statement $(ii)$ of Lemma \ref{l:c}. Let $C_1\klk C_s$ be the irreducible components of $V$. For $1\le j \le s$ we denote the dominating morphism of
affine varieties $\varphi_{|_{C_j}}: C_j\to \A^r$ by $\varphi_j$. Observe that $\deg\,\varphi=\sum_{1\le j \le s}\varphi_j$ holds. Since our ground field $\C$ is of characteristic
zero, \cite{he}, Proposition 1 (ii) implies that we may choose a non--empty Zariski open
subset $U$ of $\A^r$ satisfying the following condition:\\
For any $y\in U$, the fibers $\varphi_j^{-1}(y),\;1\le j \le s$ are of cardinality $\deg\,\varphi_j$ and form a (disjoint) partition of  $\varphi^{-1}(y)$.\\  
This implies that for any $y\in U$ 
\[
\#\,\varphi^{-1}(y)=\sum_{1\le j \le s}\#\, \varphi_j^{-1}(y)=\sum_{1\le j \le s}\deg\,\varphi_j=\deg\,\varphi.
\]
holds.
\end{prf}
\begin{observation}\label{o:b}
Statement and proof of Lemma \ref{l:c} $(i)$ do not depend on the characteristic of the ground field (in our case $\C$).
\end{observation}

Now we are going to prove Theorem \ref{th:c}.\spar
\begin{prf}
We limit our attention to the classic polar varieties. The case of the
dual polar varieties can be treated similarly (but not identically). In
order to simplify the exposition we assume that $S$ is smooth.\spar
Let us fix $1\le i \le n-p$. Without loss of generality we may suppose that 
$W_{\ul{K}^{n-p-i}(\ul{a}^{(i)})}(S)$ is non--empty. We consider the closed
subvariety $W_i$ of $\A^n\times \A^{(n-p-i+1)\times n}$ defined by the vanishing of $\Fop$ and of all $(n-i+1)$--minors of the $((n-i+1)\times n)$--matrix
\[
P_i:=\begin{bmatrix}
&J\iFop&\\
A_{1,1}&\cdots&A_{1,n}\\
\vdots&\vdots&\vdots\\
A_{n-p-i+1,1}&\cdots&A_{n-p-i+1,n}\\
\end{bmatrix}
\] 
and the canonical projection $\pi_i:=W_i\to \A^{(n-p-i+1)\times n}$ of $W_i$ to $\A^{(n-p-i+1)\times n}$. Let $D$ be an irreducible component of $W_{\ul{K}^{n-p-i}(\ul{a}^{(i)})}(S)$. Then there exists an irreducible component $C$ of $W_i$ which
contains $D\times \{a^{(i)}\}$. We claim that there exists a $((n-i)\times (n-i))$--
submatrix of $P_i$ containing $(n-i)$ entries of each of the rows number $1\klk p$ of $P_i$, such that the corresponding $(n-i)$--minor of $P_i$, say
\[
m:=\det \begin{bmatrix}
\frac{\partial F_1}{\partial X_1}&\cdots&\frac{\partial F_1}{\partial X_{n-i}}\\
\vdots&\cdots&\vdots\\
\frac{\partial F_p}{\partial X_1}&\cdots&\frac{\partial F_p}{\partial X_{n-i}}\\
A_{1,1}&\cdots&A_{1,n-i}\\
\vdots&\vdots&\vdots\\
A_{n-p-i,1}&\cdots&A_{n-p-i,n-i}\\
\end{bmatrix},
\]
does not vanish identically on $C$.\\
Otherwise $C$ would be contained in the locus defined by the vanishing of $\Fop$ and 
the determinants of all $((n-i)\times (n-i))$--submatrices of $P_i$ that contain $n-i$
entries of each of the rows number $1\klk p$ of $P_i$. This would imply that $D$ is contained in $W_{\ul{K}^{n-p-i-1}(\ul{a}^{(i+1)})}(S)$ and hence 
\[
\dim\,D\le \dim\, W_{\ul{K}^{n-p-i-1}(\ul{a}^{(i+1)})}(S) = n-p-i-1.
\]
Since by assumption the codimension of $D$ in $S$ is $i$, we have $\dim\, D=n-p-i$
and therefore a contradiction.\\
Thus we may assume without loss of generality that the
$(n-i)$--minor $m$ of the $((n-i+1)\times n)$--matrix $P_i$ does not vanish identically 
on $C$.\spar
For $n-i< j \le n$ let us denote by $M_j$ the $(n-i+1)$--minor of the $((n-i+1)\times n)$--
matrix $P_i$ which corresponds to the first $n-i$ columns and the column number $j$ of $P_i$.\\
Let $\widetilde{W}_i$ be the closed subvariety of $\A^n\times \A^{(n-p-i+1)\times n}$
where $\Fop$ and the $(n-i+1)$--minors $M_{n-i+1}\klk M_n$ vanish. Observe that any
irreducible component of $\widetilde{W}_i$ has dimension at least 
$n+(n-p-i+1)n-p-i=(n-p-i+2)n -p -i\ge 0$. 
Obviously $\widetilde{W}_i$ contains $W_i$ and therefore $C$. Thus there exists an irreducible component $\widetilde{C}$ of $\widetilde{W}_i$ that contains $C$. Observe that $m$ does not vanish identically on $\widetilde{C}$ and 
that $\widetilde{C}_m\supset C_m\neq \emptyset$ holds. From the Exchange Lemma in \cite{bank2}
we deduce now that all $(n-i+1)$--minors of the $((n-i+1)\times n)$--matrix $P_i$ vanish identically on $\widetilde{C}_m$ and hence on $\widetilde{C}$. Thus $\widetilde{C}$ is an irreducible closed subset of $W_i$ that contains the component $C$. This implies $\widetilde{C}=C$. 
Therefore $C$ is an irreducible component of the variety $\widetilde{W}_i$, whence
$\dim\,C \ge (n-p-i+2)n -p-i$.\spar
We claim now that the morphism of affine varieties 
${\pi_i}_{|_{C}}: C \to \A^{(n-i+1)\times n}$
 is dominating. Otherwise, by the Theorem on Fibers, any irreducible component of any (non--empty) fiber of ${\pi_i}_{|_{C}}$ would have a dimension strictly larger than 
\[
\dim\,C -(n-p-i+1)n \ge ((n-p-i+2)n -p-i)-(n-p-i+1)n=n-p-i.
\]
Recall that $D$ is an irreducible component of $W_{\ul{K}^{n-p-i}(\ul{a}^{(i)})}(S)$ such that $D\times \{a^{(i)}\}$ is contained in $C$. Moreover we have
$\pi^{(-1)}_i(a^{(i)})\cong W_{\ul{K}^{n-p-i}(\ul{a}^{(i)})}(S)$. Thus $D\times \{a^{(i)}\}$ is an irreducible component of $\pi^{(-1)}_i(a^{(i)})$
which is contained in $\pi^{(-1)}_i(a^{(i)})\cap C$. Since by assumption $W_{\ul{K}^{n-p-i}(\ul{a}^{(i)})}(S)$ is of pure codimension $i$ in $S$, we have $\dim\, D=n-p-i$ and therefore the ${\pi_i}_{|_C}$--fiber of $a^{(i)}\in \A^{(n-p-i+1)\times n}$, namely $\pi_i^{-1}(a^{(i)})\cap C$, contains an irreducible component 
of dimension exactly $n-p-i$, contradicting our previous conclusion that the irreducible components of any ${\pi_i}_{|_C}$--fiber are all of dimension strictly larger than $n-p-i$.\spar
Therefore the the morphism of affine varieties ${\pi_i}_{|_C}:C \to \A^{(n-i+1)\times n}$ 
 is dominating. Taking into account the estimate  $\dim\,C\le (n-p-i+2)n-p-i$
 and that the ${\pi_i}_{|_C}$--fiber of $a^{(i)}$ contains an irreducible component of dimension $n-p-i$, we deduce now from the Theorem on Fibers that $\dim\,C= (n-p-i+2)n-p-i$ holds.
 \spar 
 Let $V_i$  be the union of all $((n-i+2)n-p-i)$--dimensional components $C$
 of $W_i$ such that the morphism of affine varieties  ${\pi_i}_{|_C}:C \to \A^{(n-p-i+1)\times n}$ is dominating. Let us denote the restriction map  ${\pi_i}_{|_{V_i}}:V_i \to \A^{(n-p-i+1)\times n}$ by $\psi_i:V_i\to \A^{(n-p-i+1)\times n}$.
\spar
  
Then $V_i$ is an equidimensional affine variety of dimension $(n-p-i+2)n-p-i$ and $\psi_i$ is an morphism of affine varieties such that the restriction of $\psi_i$ to any irreducible component of $V_i$ is dominating. Observe that our previous argumentation implies 
 \[
 \psi_i^{-1}(a^{(i)})\cong W_{\ul{K}^{n-p-i}(\ul{a}^{(i)})}(S)\;\;
 \text{and}\;\;
 \psi_i^{-1}(b^{(i)})\cong W_{\ul{K}^{n-p-i}(\ul{b}^{(i)})}(S)\;\;
 \text{for any}\;\; b\in U. 
 \]
Observe now that there exists a finite set $\mathcal{M}$ of rational $((n-p-i)\times n)$--matrices of maximal rank $n-p-i$ satisfying the following conditions:
\begin{itemize}
\item [-] Each $M\in \mathcal{M}$ represents a generic Noether position (see \cite{he}, Lemma 1) of the polar variety $W_{\ul{K}^{n-p-i}(\ul{a}^{(i)})}(S)$.
\item [-] For each $b\in U$ there exists a $((n-p-i)\times n)$ matrix $M\in \mathcal{M}$
representing a generic Noether position of the affine polar variety $W_{\ul{K}^{n-p-i}(\ul{b}^{(i)})}(S)$ (here we use the assumption that $W_{\ul{K}^{n-p-i}(\ul{a}^{(i)})}(S)\neq \emptyset$ implies  $W_{\ul{K}^{n-p-i}(\ul{b}^{(i)})}(S)\neq \emptyset$).
\end{itemize}
Therefore, restricting, if necessary, the non--empty Zariski open set $U$, we may suppose without loss of generality that there exists a full--rank matrix $M\in \Q^{(n-p-i)\times n}$
such that $M$ represents a generic Noether position of $W_{\ul{K}^{n-p-i}(\ul{a}^{(i)})}(S)$ and of any $W_{\ul{K}^{n-p-i}(\ul{b}^{(i)})}(S)$, where $b$ belongs to $U$. The $((n-p-i)\times n)$--matrix $M$ and $\psi_i$ induce now a morphism of $((n-p-i+2)n -p-i)$--dimensional affine varieties $\varphi_i:V_i\to \A^{n-p-i}\times \A^{(n-p-i+1)\times n}$ which associates to any point $(x,c)\in V_i\subset \A^n\times \A^{(n-p-i+1)\times n}$ the value $\varphi_i(x,c):=(Mx,c)$.\spar 

Let $C$ be an arbitrary irreducible component of $V_i$. Recall that we have $\dim\,C=(n-p-i+2)n-p-i$ and the morphism of affine varieties ${\psi_i}_{|_C}:C\to \A^{(n-p-i+1)\times n}$ is dominating. Therefore, by the Theorem on Fibers, there exists  $b\in U$ such that $b^{(i)}$ belongs to $\psi_i(C)$ and that the ${\psi_i}_{|_C}$--fiber of
$b^{(i)}$, namely $\psi_i^{-1}(b^{(i)}\cap C$,
is non--empty and equidimensional of dimension $(n-p-i+2)n-p-i-(n-p-i+1)n=n-p-i$.
Observe that $b\in U$ implies that $W_{\ul{K}^{n-p-i}(\ul{b}^{(i)})}$ is non--empty and of pure codimension $i$ in $S$. Therefore any irreducible 
component $W_{\ul{K}^{n-p-i}(\ul{b}^{(i)})}$ is of dimension $n-p-i$.
Hence any irreducible component of $\psi_i^{-1}(b^{(i)})\cap C$ is isomorphic to an irreducible component of $W_{\ul{K}^{n-p-i}(\ul{b}^{(i)})}(S)$. This implies that the rational $((n-p-i)\times n)$--matrix $M$ represents 
also a generic Noether position of the ${\psi_i}_{|_C}$--fiber of
$b^{(i)}$, namely $\psi_i^{-1}(b^{(i)})\cap C$.\spar

Fix now any point $y\in \A^{n-p-i}$ and observe that the ${\varphi_i}_{|_C}$--fiber
$(y,b^{(i)})$ is isomorphic to a non--empty subset of the zero--dimensional variety 
\[
\{x\in W_{\ul{K}^{n-p-i}(\ul{b}^{(i)})}(S)\;|\; Mx=y\}.
\]
Therefore we may conclude that the morphism of affine varieties 
\[
{\varphi_i}_{|_C}:C\to \A^{n-p-i}\times \A^{(n-p-i+1)\times n}
\]
is dominating. Since $C$ was chosen as an arbitrary irreducible component of $V_i$,
we infer that $\varphi_i:V_i\to \A^{n-p-i}\times \A^{(n-p-i+1)\times n}$ satisfies the
requirements of Lemma \ref{l:c}.\spar

Let us now fix $b\in U$ and $y\in \A^{n-p-i}$ such that $(y, b^{(i)})$ represents a generic choice in the ambient space $ \A^{n-p-i}\times \A^{(n-p-i+1)\times n}$. Observe that ${\varphi_i}^{-1}(y,a^{(i)})$ is isomorphic to
\[
\{x\in W_{\ul{K}^{n-p-i}(\ul{a}^{(i)})}(S)\;|\; Mx=y\}.
\]
Since $M$ represents a generic Noether position of the affine variety $W_{\ul{K}^{n-p-i}(\ul{a}^{(i)})}(S)$, we conclude that $\#\,\varphi_i^{-1}(y,a^{(i)})=\deg\,W_{\ul{K}^{n-p-i}(\ul{a}^{(i)})}(S)$ holds. In particular,\\
$\#\,\varphi_i^{-1}(y,\ul{a}^{(i)})$ is finite.
Similarly, one deduces $\#\,\varphi_i^{-1}(y,b^{(i)})=\deg\,W_{\ul{K}^{n-p-i}(\ul{b}^{(i)})}(S)$.\spar

Lemma \ref{l:c} implies now 
\[
\deg\,W_{\ul{K}^{n-p-i}(\ul{b}^{(i)})}(S)=\#\,\varphi_i^{-1}(y,b^{(i)})=\deg\, \varphi_i
\]
and
\[
\deg\,W_{\ul{K}^{n-p-i}(\ul{a}^{(i)})}(S)=\#\,\varphi_i^{-1}(y,a^{(i)})\le \deg\, \varphi_i=\deg\,W_{\ul{K}^{n-p-i}(\ul{b}^{(i)})}(S).
\]
This proves Theorem \ref{th:c}.
\end{prf}

The following remark is now at order.
Theorem \ref{th:c} is not aimed to replace a general intersection theory formulated in terms of rational equivalence classes of cycles of the projective space $\P^n$, it is rather complementary to such a theory. 
As said before, its main purpose is to allow, in terms of the degrees of generic polar varieties, complexity estimates for real point finding procedures  
which are based on the consideration of meagerly generic polar varieties. The alluded
interplay between a general intersection theory and the statement of  Theorem \ref{th:c}
may become more clear by the following considerations.\spar

Under the assumption that the homogenizations of the polynomials $F_1\klk F_p$
define a smooth codimension $p$ subvariety $\widetilde{S}$ of $\P^n$, the projective 
closures of the affine polar varieties $W_{\ul{K}(\ul{a}^{(i)})}(S)$, $1\le i \le n-p$,
form the corresponding polar varieties of $\widetilde{S}$, when $a\in \A^{(n-p)\times n}$ is generic. Moreover, $\widetilde{S}$ is the projective closure of $S$.
The degrees of the generic (classic) polar varieties of $\widetilde{S}$
may be expressed in terms of the degrees of the Chern classes of $\widetilde{S}$
(see \cite{fu}, Example 14.3.3). Further, since $\widetilde{S}$ is a projective smooth 
complete intersection variety, the total 
Chern class of $\widetilde{S}$ may be characterized in terms of the degrees 
of $F_1\klk F_p$ and the first Chern class $c_1$ of the normal bundle of
$\widetilde{S}$ in $\P^n$ (see \cite{hi}, Theorem 4.8.1, Section 22.1 and 
\cite{nav}, Theorem 1). This implies an upper bound for the geometric degrees of the {\em affine} polar varieties $W_{\ul{K}(\ul{a}^{(i)})}(S)$ in terms of the degrees of the polynomials $F_1\klk F_p$ and the class $c_1$.\spar

Of course our previous assumption on $F_1\klk F_p$ is very restrictive. Nevertheless,
in this situation we are able to illustrate how results like Theorem \ref{th:c} (or \cite{he}, Proposition 1) interact with facts from a general intersection theory.

\spar
Let us now turn back to the discussion of the, up to now informal, concept of meagerly generic polar varieties. We are going to give to this concept a precise mathematical shape
and  to discuss it by concrete examples.\spar

\begin{definition}\label{d:1}
Let $1\le i \le n-p$ and let $m$ be a non--zero polynomial of $\C[X]$, let $E$ be an irreducible closed subvariety of $\A^{(n-p-i+1)\times n}$ and let $O$ be a non--empty
Zariski open subset of $E$. Suppose that the following conditions are satisfied:
\begin{itemize}
\item[(i)] each $b\in O$ is a complex $((n-p)\times n)$--matrix of maximal rank $n-p-i+1$,
\item[(ii)] for each $b\in O$ the affine variety $W_{\ul{K}(\ul{b})}(S)_{m}$ (respectively $W_{\ol{K}(\ol{b})}(S)_{m}$) is empty or of pure
codimension $i$ in $S$.
\end{itemize}
Then we call the algebraic family $\{W_{\ul{K}(\ul{b})}(S)_{m}\}_{b\in O}$ (respectively $\{W_{\ol{K}(\ol{b})}(S)_{m}\}_{b\in O}$) of Zariski open subsets of polar varieties of $S$ {\em meagerly generic}.\spar

In the same vein, we call for $b\in O$ the affine variety $W_{\ul{K}(\ul{b})}(S)_{m}$ (respectively $W_{\ol{K}(\ol{b})}(S)_{m}$) meagerly generic.
\end{definition}

Typically, the irreducible affine variety $E$ contains a Zariski dense set of rational points and for any
$b\in O$ the non--emptiness of $W_{\ul{K}(\ul{b})}(S)_{m}$ (respectively $W_{\ol{K}(\ol{b})}(S)_{m}$) implies the non--emptiness of the open locus 
of the corresponding generic polar variety defined by the non-vanishing of $m$. From Theorem \ref{th:c}
and its proof we deduce that in this case the sum of the degrees of the irreducible components of $W_{\ul{K}(\ul{b})}(S)$ (respectively $W_{\ol{K}(\ol{b})}(S)$), where $m$ does not vanish identically, is bounded by the sum of the degrees 
of the components of the corresponding generic polar variety which satisfy the same condition.\spar

We paraphrase this consideration as follows.

\begin{observation}\label{o:c}
For full rank matrices $\A^{(n-p)\times n}$ and $1\le i \le n-p$ the degrees 
of (suitable open loci) of meagerly generic polar varieties $W_{\ul{K}(\ul{a}^{(i)})}(S)$ and $W_{\ol{K}(\ol{a}^{(i)})}(S)$ attain their maximum when $a$ is generic.
\end{observation}

We are now going to discuss two examples of geometrically relevant algebraic families of Zariski open subsets of polar varieties of $S$. It will turn out that these sets are empty or smooth affine subvarieties of $S$.\bigskip

\textbf{Example 1}\spar

We are going to adapt the argumentation used in \cite{bank2}, Section 2.3
to the terminology of meagerly generic polar varieties.\spar
Let $m \in \Qxon$ denote the $(p-1)$--minor of the Jacobian 
$J(F_1\kpk F_p)$ given by the first $(p-1)$ rows and columns, i.e., let 
\[
m := \det \;
\left[ \displaystyle\displaystyle\frac{\partial F_k}{\partial X_l}\right]_
{{1 \le k \le p-1} \atop {1 \le l \le p-1}}. 
\]
Moreover, for $p\le r \le n,\, p\le t < n$ let $Z_{r,t}$
be a new indeterminate.
Using the following regular $((n-p+1) \times (n-p+1))$--parameter matrix
\[
Z:= \begin{bmatrix}
1 & 0 &\cdots & 0 &\cdots&&\cdots & 0\\
Z_{p+1,p} & 1 &&&&&&\\
\vdots& \vdots &  \ddots & &\mbox{\Huge O}   && &\vdots \\
Z_{p+i-1,p} & Z_{p+i-1,p+1} &\cdots &  1   &&&&\\
Z_{p+i,p}  &  Z_{p+i,p+1} &\cdots & Z_{p+i, p+i-1}  &  1 &&&\\
\vdots & \vdots &  & \vdots & \vdots &\ddots && 0 \\
Z_{n,p} & Z_{n,p+1} & \cdots & Z_{n,p+i-1} & Z_{n,p+i} & Z_{n,p+i+1}&
\cdots  &1
\end{bmatrix},
\]
we are going to introduce an $(n \times n)$--coordinate
transformation matrix $A:=A(Z)$
which will represent the key for our ongoing construction of a
meagerly generic family of Zariski open subsets of classic polar varieties of
$S$.\smallskip

Let us fix an index $1 \le i \le n-p$. According
to our choice of $i$,
the matrix $Z$ may be subdivided into submatrices as follows:
\[
Z =\begin{bmatrix}
Z_1^{(i)} & O_{i,n-p-i+1} \\  Z^{(i)} & Z_2^{(i)}
\end{bmatrix}.
\]
Here the matrix $Z^{(i)}$ is defined as
\[
Z^{(i)} := \begin{bmatrix}
Z_{p+i,p} & \cdots & Z_{p+i,p+i-1} \\
\vdots & \cdots & \vdots \\
Z_{n,p} & \cdots &Z_{n,p+i-1}
\end{bmatrix},
\]
and $Z^{(i)}_1$ and $Z_2^{(i)}$ denote the quadratic lower triangular
matrices bordering $Z^{(i)}$ in $Z$.
Let
\[
A:=\begin{bmatrix} I_{p-1} & O_{p-1,i} &
O_{p-1,n-p-i+1}\\ O_{i,p-1} & Z_1^{(i)} & O_{i,n-p-i+1}\\
O_{n-p-i+1,p-1} & Z^{(i)} & Z_2^{(i)} \end{bmatrix}.
\]
Here the submatrix $I_{p-1}$ is the $((p-1) \times (p-1))$--identity matrix and
$Z^{(i)},
Z_1^{(i)},$ and $Z_2^{(i)}$ are the submatrices of the parameter matrix $Z$
introduced before. Thus, $A$ is a regular, parameter dependent
$(n \times n)$--coordinate transformation matrix.
\spar

Like the matrix $Z$, the matrix $A(Z)$ contains
\[
s:=\frac{(n-p) (n-p+1)}{2}
\]
parameters $Z_{r,t}$ which we may specialize into
any point $z$ of the affine space $\A^s$. For such a point $z \in \A^s$ we
denote
the corresponding specialized matrices by $Z_1^{(i)}(z), Z_2^{(i)}(z), 
Z^{(i)}(z)$ and $A(z)$.\spar

Let $B_i:= B_i(Z)$ be the $((n-p-i+1) \times n)$--matrix consisting of the
rows
number $p+i \klk n$ of the inverse matrix $A(Z)^{-1}$ of $A(Z)$. From the
particular triangular form of the matrix $A(Z)$ we deduce that the entries
of $B_i(Z)$
are polynomials in the indeterminates $Z_{r,t}$, $p\le r \le n,\,p \le t < r$.
Moreover, we have the matrix identity
\[B_iA=
\begin{bmatrix} O_{n-p-i+1,p+i-1} & I_{n-p-i+1} \end{bmatrix}.
\]
\spar

Finally, let $E$ be the Zariski closure of $B_i(\A^s)$ in $\A^{(n-p-i+1) \times n}$.
Then $E$ is an irreducible closed subvariety of $\A^{(n-p-i+1) \times n}$
consisting of complex $((n-p-i+1) \times n)$--matrices of maximal rank
$n-p-i+1$.
\spar

According to the structure of the coordinate transformation matrix
$A=A(Z)$ we
subdivide the Jacobian $J(F_1 \klk F_p)$ into three
submatrices
\[
J(F_1 \klk F_p)=\begin{bmatrix} U & V^{(i)} & W^{(i)} \end{bmatrix}, \]
with
\[
U:=\left[ \displaystyle\displaystyle\frac{\partial F_k}{\partial
X_l}\right]_
{{1 \le k \le p} \atop {1 \le l \le p-1}},\;
V^{(i)}:=\left[ \displaystyle\displaystyle\frac{\partial F_k}{\partial
X_l}\right]_
{{1 \le k \le p} \atop {p \le l \le p+i-1}},\;
W^{(i)}:=\left[ \displaystyle\displaystyle\frac{\partial F_k}{\partial
X_l}\right]_
{{1 \le k \le p} \atop {p+i \le l \le n}}.
\]
We have then the following matrix identity:
\[
J(F_1 \klk F_p)\;A(Z) =
\begin{bmatrix} U\; & V^{(i)}Z_1^{(i)}+W^{(i)}Z^{(i)} &\; W^{(i)}Z^{(i)}_2
\end{bmatrix}.
\]

For $p\le j \le p+i-1$ let us denote by $\widetilde{M}_j :=
\widetilde{M}_j(X,Z)$ the $p$--minor of the polynomial $(p \times
n)$--matrix $J(F_1 \klk F_p)\;A(Z)$ that corresponds to the
columns number $1 \klk p-1,j$.
\spar

From \cite{bank2}, Section 2.3 and Lemma 1 we deduce that there exists a
non--empty Zariski open subset $\widetilde{O}$ of $\A^s$ such that any
point $z \in \widetilde{O}$ satisfies the following two conditions:
\begin{itemize}
\item[$(i)$] the polynomial equations 
\[F_1(X)=0 \klk
F_p(X)=0,\widetilde{M}_p(X,z)=0 \klk \widetilde{M}_{p+i-1}(X,z)=0\]
intersect transversally at any of their common solutions in $\A^n_m$,
\item[$(ii)$] for any point $x \in \A^n_m$ with 
\[
F_1(x)=0 \klk  F_p(x)=0,\widetilde{M}_p(x,z)=0 \klk \widetilde{M}_{p+i-1}(x,z)=0
\]
all $p$--minors of the complex $(p \times (p+i-1))$--matrix 
\[
\left[
\;U(x)\;\;V^{(i)}(x)Z_1^{(i)}(z)+W^{(i)}(x)Z^{(i)}(z)\right]
\] 
vanish, i.e., this matrix has rank at most $p-1$.
\end{itemize}

Let $O:= B_i(\widetilde{O})$ and observe that $O$ is a non--empty Zariski
open subset of the irreducible affine variety $E$.
\spar

Let us fix for the moment an arbitrary point $z \in \widetilde{O}$. Then
the complex $((n-p-i+1) \times n)$--matrix $B_i(z)$ is of maximal rank
$n-p-i+1$.\spar

We consider now the Zariski open subset of $W_{\ul{K}(\ul{B_i(z)})}(S)_m$ of
the classic polar variety\\ $W_{\ul{K}(\ul{B_i(z)})}(S)$. By
definition the affine variety $W_{\ul{K}(\ul{B_i(z)})}(S)_m$
consists of the points of $S_m$ where all $(n-i+1)$--minors of the
polynomial $((n-i+1) \times n)$--matrix
\[
N_i:=\begin{bmatrix}J\iFop\\B_i(z)\end{bmatrix}
\]
vanish.\spar

Let $x$ be an arbitrary element of $S_m$. Then all $(n-i+1)$--minors of $N_i$
vanish at $x$ if and only if $N_i$ has rank at most $n-i$ at $x$. This
is equivalent to the condition that 
\[
N_i\,A(z)=\begin{bmatrix}J\iFop\,A(z)\\B_i(z)\,A(z)\end{bmatrix}
=\begin{bmatrix}
U&V^{(i)}Z_1^{(i)}(z)+W^{(i)}Z^{(i)}(z)&W^{(i)}Z^{(i)}_2(z)\\&O_{n-p-i+1,p+i-1}&I_{n-p-i+1}
\end{bmatrix}
\]
has rank at most $n-i$ at $x$. Hence
$W_{\ul{K}(\ul{B_i(z)})}(S)_m$ consists of the points of $S_m$
where the $(p \times (p+i-1))$--matrix \[ \left[
\;U\;\;V^{(i)} Z_1^{(i)}(z)+W^{(i)}Z^{(i)}(z)\right]\] has rank at most $p-1$.
\spar

From conditions $(i)$ and $(ii)$ above we deduce that the affine
variety $W_{\ul{K}(\ul{B_i(z)})}(S)_m$ is either empty or a smooth
subvariety of $S$ of pure codimension $i$.\spar

Since $z \in \widetilde{O}$ was chosen arbitrarily, we have shown the
following restatement of \cite{bank2}, Theorem 1:
\spar
\begin{proposition}\label{p:c}
For each point $b \in O$ the Zariski open subset
$W_{\ul{K}(\ul{b})}(S)_m$  of the polar variety
$W_{\ul{K}(\ul{b})}(S)$ is either empty or a smooth affine
subvariety of $S$ of pure codimension $i$. Thus
$\{W_{\ul{K}(\ul{b})}(S)_m\}_{b \in O}$  forms a meagerly generic
algebraic family of empty or smooth Zariski open subsets of polar varieties of $S$.
\end{proposition}
Observe that the matrix $A(Z)$ (and hence the matrix $A(Z)^{-1}$) does not depend on the index $1\le i \le n-p$. Moreover,  $A(Z)^{-1}$ has the same triangular shape as $A(Z)$. Here the entries $Z_{r,t}\,,\,p\le r \le n,\, p\le t < r$ have to be replaced by suitable polynomials in $Z$ which are algebraically independent over $\Q$.\spar

Consequently
\[
B_{n-p}(Z)\klk B_1(Z)
\]
form a nested sequence of matrices, each contained in the other. Therefore we obtain for
$z\in \widetilde{O}$ a descending chain of Zariski open susets of polar varieties
\[
W_{\ul{K}(\ul{B_1}(z))}(S)\,_m \supset \cdots \supset W_{\ul{K}(\ul{B_{n-p}}(z))}(S)\,_m
\]
which are empty or smooth codimension one subvarieties one of the other.\spar

We discuss now Example 1 in view of algorithmic applications. In order to do this let us start
with the following considerations.\spar

Let $Y^*=(Y_1\klk Y_{p-1})$ be a $(p-1)$--tupel of new indeterminates . We fix for the moment indices $1\le h_1 < \cdots < h_p\le n $ and make a sub--selection of $p-1$ of them,
say $h'_1<\cdots < h'_{p-1}$. For the sake of notational simplicity suppose $h_1=1\klk h_p=p$ and $h'_1=1\klk h'_{p-1}=p-1$.\spar

Let
\[
U^*:=\left[\frac{\partial F_k}{\partial X_l}\right]_{1\le k,l \le p},\;\;V^*:=\left[\frac{\partial F_k}{\partial X_l}\right]_{1\le k \le p \atop{p+1\le l \le n}}, M^*:=\det U^*
\]
and let $N^*:=N^*(X,Y^*)$ be the polynomial $((p+1)\times p)$--matrix
\[
N^*:=\begin{bmatrix}
U^*\\Y_1\cdots Y_{p-1}\;1
\end{bmatrix}.
\]
There exists a polynomial $(p\times p)$--matrix $C^*:=C^*(Y^*)$ with $\det\,C^*=1$ satisfying the condition 
\[
N^*\, C^*=\begin{bmatrix}
U^*\, C^*\\
0\cdots 0\;1
\end{bmatrix}.
\]
In particular we have $M^*= det(U^*C^*)$.\spar

To our sub--selection of indices corresponds the $(p-1)$--minor $m^*:=m^*(X,Y^*)$
determined by the rows and columns number $1\klk p-1$ of $U^*\,C^*$. For the moment let us fix an index $1\le i \le n-p$.
According to our choice of indices and the sub--selection made from them,
let us denote for $p\le j \le p+i-1$ by $M_j^*=M_j^*(X,Y^*,Z)$ the
$p$--minor of the polynomial $(p\times n)$--matrix
$\left[U^*\,C^*\;\;V^*\right]\,A(Z)$ that corresponds to the columns
number $1\klk p-1,j$.
\spar

We consider now the polynomial map
\[
\Phi^*: \A^n_{M^*}\times \A^{p-1}\times \A^s \to \A^p\times \A^i 
\]
defined for $(x,y^*,z)\in \A^n_{M^*}\times \A^{p-1}\times \A^s$ by
\[
\Phi^*(x,y^*,z):=(F_1(x)\klk F_p(x), M^*_p(x,y^*,z)\klk M^*_{p+i-1}(x,y^*,z)).
\]
As in \cite{bank2}, Section 2.3 we may now argue that $0\in \A^{p-1}\times \A^i$
is a regular value of the polynomial map $\Phi^*$. From the Weak Transversality Theorem
of Thom--Sard we then deduce that there exists a non--empty Zariski open subset $O^*$ of
$\A^{p-1}\times \A^i$ such that for any point $(y^*,z)$ of $O^*$ the polynomial equations 
\[
F_1(X)=0\klk F_p(X)=0, M^*_p(X,y^*,z)=0\klk M^*_{p+i-1}(X,y^*,z)=0
\]
intersect transversally at any of their common solutions in $\A^n_{M^*}$.
\spar

Let $(y^*, z)$ be an arbitrary point of $O^*$.\spar

Notice that for any point $x\in S_{M^*}$ there exists a $(p-1)$--minor of $U^*\,C^*(y^*)$ that does not vanish at $x$. Let us suppose $m^*(x,y^*)\neq 0$. Now we may conclude as in the proof of Proposition \ref{p:c} that the assumption
\[
M^*_p(x,y^*,z)=0\klk M^*_{p+i-1}(x,y^*,z)=0
\]
implies that the polynomial $((n-i+1)\times n)$--matrix 
\[
\begin{bmatrix}
U^*\,C^*(y^*)\;\;V^*\\
B_{i}(z)
\end{bmatrix}
\]
has rank at most $n-i$ at $x$.\spar

Observe that the row number $n-p-i+1$ of $B_i(Z)$ has the form\\ $(0\klk 0,  
B_{n,p}(Z)\klk B_{n,n-1}(Z),1)$, where $B_{n,p}(Z)\klk B_{n,n-1}(Z)$ are suitable polynomials in $Z$ which are algebraically independent over $\Q$. Without loss of generality we may suppose that any point $(y^*,z)$ of $O^*$ satisfies the condition $B_{n,p}(z)\neq 0$. Thus there exists a non--empty Zariski open subset $Q^*$ of $\A^n$ such that any point $b=(b_1\klk b_n)$ of $Q^*$ has non--zero entries and satisfies the following 
condition:\spar
There exists a point $(y^*,z)$ of $O^*$ with $y^*=\frac{1}{b_p}(b_1\klk b_{p-1})$
and 
\[
\frac{1}{B_{n,p}(z)}(B_{n,p+1}(z)\klk B_{n,n-1}(z), 1)=\frac{1}{b_p}(b_{p+1}\klk b_n).
\]
All these constructions and arguments depend on a fixed choice of indices $1\le h_1<\cdots <h_p\le n$, namely $h_1:=1\klk h_p:=p$ and a sub--selection of them, namely $h'_1:=1\klk h'_{p-1}:=p-1$.\spar

Intersecting now for all possible choices of indices and sub--selections of them the resulting subsets of $\A^n$ corresponding to $Q^*$ in case $h_1=1\klk h_p=p,\; h'_1=1\klk h'_{p-1}=p-1$, we obtain a non--empty Zariski open subset $Q$ of $\A^n$ which has the same properties
as $Q^*$ with respect to any choice of indices $1\le h_1 <\cdots < h_p\le n$ and any sub--selection $h'_1<\cdots < h'_{p-1}$ of them.\spar

Let us now suppose that  $S_{\R}$ is smooth and compact and let $b=(b_1\klk b_n)$ be a fixed point of $Q\cap \Q^n$. Of course, such a point $b$ exists and can be obtained by a random choice.\spar

From \cite{bank2}, Theorem 2 we deduce that $W_{\ul{K}(\ul{b})}(S_{\R})$
is not empty.  Since $b$ belongs to $Q$ there exists a point $(y^*,z)\in Q^*$ with 
\[
y^*=\frac{1}{b_p}(b_1\klk b_{p-1})\;\;\text{and}\;\;\frac{1}{B_{n,p}(z)}(B_{n,p+1}(z)\klk B_{n,n-1}(z), 1)=\frac{1}{b_p}(b_{p+1}\klk b_n).
\]
We may assume  $z\in \Q^s$. Since $S_{\R}$ is smooth and $W_{\ul{K}(\ul{b})}(S_{\R})$ is not empty, we may suppose without loss of generality 
\[
W_{\ul{K}(\ul{b})}(S)_{M^*(X)m^*(X,y^*)}\neq \emptyset.
\]
For $1\le i \le n-p$ let $W_i^*$ be the affine subvariety of it defined by the vanishing of all $(n-i+1)$--minors of the polynomial $((n-i+1)\times n)$--matrix $\begin{bmatrix}
U^*\,C^*(y^*)\;\;V^*\\
B_{i}(z)
\end{bmatrix}$. As before we obtain a descending chain of affine varieties, namely 
$W^*_1\supset \cdots \supset W^*_{n-p}$ which are empty or smooth codimension one
subvarieties one of the other. In order to rule out the possibility of emptiness we are going to show $W^*_{n-p}\neq \emptyset$.\spar

For this purpose we consider an arbitrary element $x\in S_{M^*(X)m^*(X,y^*)}$. Then all $(p+1)$--minors of the polynomial $((p+1)\times n)$--matrix $\begin{bmatrix}
J\iFop\\b \end{bmatrix}$ vanish at $x$ if and only if the same is true for   
$\begin{bmatrix}
J\iFop\\\frac{1}{b_p}b \end{bmatrix}$ and this is equivalent to the condition that the polynomial $((p+1)\times n)$--matrix
\[
\begin{bmatrix}
U^*\,C^*(y^*)\;\;V^*\\
0\cdots 0\;1\;\frac{1}{b_p}b_{p+1}\cdots \frac{1}{b_p} b_n
\end{bmatrix},
\]
and hence also
\[
\begin{bmatrix}
U^*\,C^*(y^*)\;\;V^*\\
0\cdots 0\;b_p\;b_{p+1}\cdots  b_n
\end{bmatrix}=\begin{bmatrix}
U^*\,C^*(y^*)\;\;V^*\\
0\cdots 0\;B_{n,p}(z)\cdots  B_{n,n-1}(z),\,1
\end{bmatrix},
\]
have rank at most $p$ at $x$. This implies 
\[
W_{\ul{K}(\ul{b})}(S)_{M^*(X)m^*(X,y^*)}=W^*_{n-p}\,,
\quad \text{whence}\;\; W^*_{n-p}\neq \emptyset.
\]
In order to retrieve finitely many real algebraic sample points for the connected components of the real algebraic variety $S_{\R}$ we may now proceed 
by applying the general algorithm described in \cite{bank2, bank3} or \cite{bank4} as follows:
For each choice of indices $1\le h_1<\cdots < h_p\le n$ and any
sub--selection $h'_1<\cdots < h_{p-1}$ of them we generate the equations and the inequation  
that define the descending chain of affine varieties corresponding to 
$W^*_1\supset \cdots \supset W^*_{n-p}$ in order to find the (finitely many) real points contained in the last one.\spar

The set of real points obtained in this way is by virtue of \cite{bank2}, Theorem 2 a set of sample points for the connected components of $S_{\R}$. For details of the algorithm we refer to \cite{bank2}.\spar

We are now going to explain how the affine varieties $W^*_i\,,1\le i \le n-p$ may be interpreted as Zariski open subsets of polar varieties of suitable complete intersection varieties.\spar  

Suppose we have already fixed indices $1\le h_1 <\cdots < h_p \le n$ and have made a sub--selection $h'_1<\cdots <h_{p-1}$. For the sake of notational simplicity we suppose again 
$h_1=1\klk h_p=p,\,h'_1=1\klk h'_{p-1}=p-1$. Furthermore, we suppose that we have chosen 
a point $(y^*,z)$ of $O^*$. Let $\Omega=(\Omega_1\klk \Omega_n)$ be an $n$--tupel of new indeterminates and let us consider the coordinate transformation matrix 
\[D:=\begin{bmatrix}
C^*(y^*)&O_{p,n-p}\\O_{n-p,p}&I_{n-p}\end{bmatrix}\] and the polynomials
\[
G_1:=G_1(\Omega):=F_1(\Omega\,D^T)\klk G_p:=G_p(\Omega):=F_p(\Omega\,D^T),
\]
where $D^T$ denotes the transposed matrix of $D$. Let $G:=(G_1\klk G_p)$, $S_G:=\{G_1=0\klk G_p=0\}$ and $M_G$ the $(p\times p)$--minor of $J(G_1\klk G_p)$ which corresponds to the columns number $1\klk p$.\spar

From the identity 
\[
J(G)(X(D^{-1})^T)=J(F)D=\left[U^*\,C^*(y^*)\;\,V^*\right]
\]
we deduce now that for any point $x\in \A^n$ and any index $1\le i \le n-p$ the conditions
\[
x(D^{-1})^T\in W_{\ul{K}(\ul{B_i(z)})}(S_G)_{M_G}
\;\;\text{and}\;\; x\in W^*_i\]
are equivalent. Therefore $W_i^*$ is isomorphic to the polar variety $W_{\ul{K}(\ul{B_i(z)})}(S_G)$.\spar

\textbf{Example 2}\spar

We are going to consider simultaneously two meagerly generic algebraic families of Zariski open subsets of polar varieties of $S$, one in the classic and the other in the dual case.\spar

Let $m$ be the $(p-1)$--minor of $J\iFop$ introduced in Example 1, namely 
$m:=\det\,\left[\frac{\partial F_k}{\partial X_l}\right]_{1\le k,l \le p-1}$. We fix
now an index $1\le i \le n-p$ and an $(n-i)$--tupel of rational numbers $\gamma=(\gamma_1\klk \gamma_{n-i})$ with $\gamma_{n-i}\neq 0$.\spar

Let $Z=(Z_{n-i+1}\klk Z_n)$ be an $i$--tupel of new indeterminates and let $B=B_{(i,\gamma)}$ be the polynomial $((n-p-i+1)\times n)$--matrix defined by
\[
B:=\begin{bmatrix}
O_{n-p-i,\,p-1}&I_{n-p-i}&O_{n-p-i,\,i+1}\\
\gamma_1\ldots \gamma_{p-1}&\gamma_p\ldots \gamma_{n-i-1}&\gamma_{n-i}\;\;
Z_{n-i+1}\ldots Z_{n}
\end{bmatrix}.
\]
By $\ul{B}:=\ul{B}_{(i,\gamma)}$ and $\ol{B}:=\ol{B}_{(i,\gamma)}$ we denote
the polynomial $((n-p-i+1)\times (n+1))$--matrices obtained by adding to $B$ as column number zero the transposed $(n-p-i+1)$--tupel $(0\klk 0,0)^T$ and $(0\klk 0, 1)^T$, respectively. So we have $\ul{B}_{\ast}=\ol{B}_{\ast}=B$.\spar

Let $E:=B(\A^i)$. Then $E$ is a closed irreducible subvariety of $\A^{(n-p-i+1)\times n}$ which is isomorphic to $\A^i$.\spar

From $\gamma_{n-i}\neq0$ we deduce that for any instance $z\in \A^i$ with $z=(z_{n-i+1}\klk z_n)$ the complex $((n-p-i+1)\times n)$--matrix $B_i(z)=\ul{B}(z)_{\ast}=\ol{B}(z)_{\ast}$
has maximal rank $n-p-i+1$ and this gives rise to two polar varieties of $S$, one classic and the other dual, namely $W_{\ul{K}(\ul{B(z)})}(S)$ and $W_{\ol{K}(\ol{B(z)})}(S)$.\spar

We shall limit our attention to the algebraic family of Zariski open subsets of dual polar varieties 
$\{W_{\ol{K}(\ol{B(z)})}(S)_{m}\}_{z\in \A^i}$. The case of the algebraic family
$\{W_{\ul{K}(\ul{B(z)})}(S)_{m}\}_{z\in \A^i}$ is treated similarly.\spar

We consider now the polynomial $((n-i+1)\times n)$--matrix $N:=N(i,\gamma)$ defined by
\[
N\!:=\!\!\scriptsize\begin{bmatrix}
J\iFop\\\begin{matrix}\\
O_{n-p-i,\,p-1}&I_{n-p-i}&O_{n-p-i,\,i+1}\\
&&\\
&&\\
\gamma_1-X_1\ldots \gamma_{p-1}-X_{p-1}&\gamma_p-X_p\ldots \gamma_{n-i-1}-X_{n-i+1}&\gamma_{n-i}-X_{n-i}\;\;
Z_{n-i+1}-X_{n-i+1}\ldots Z_{n}-X_n
\end{matrix}\end{bmatrix}
\]
Observe that the $(n-i)$--minor of $N$, which corresponds to the rows and columns $1\klk n-i$, has value $m$. For $n-i+1\le j \le n$ let us denote by $\widetilde{M}_j:=\widetilde{M}_j(X,Z)$
the $(n-i+1)$--minor of $N$ that corresponds to the columns number $1\klk n-i$ and $j$. One verifies immediately for $n-i+1\le j' \le n$ the identities $\frac{\partial \widetilde{M}_j}{\partial Z_{j'}}=0$ in case $j'\neq j $ and $\frac{\partial \widetilde{M}_j}{\partial Z_{j'}}=m$ in case $j'= j$. \spar

Let us now consider the polynomial map $\Phi: \A_m^n\times \A^i \to \A^{p+i}$ defined for $(x,z)\in \A_m^n\times \A^i$ by $\Phi(x,z):=(F_1(x)\klk F_p(x), \widetilde{M}_{n-i+1}(x,z)\klk\widetilde{M}_n(x,z))$. The Jacobian $J(\Phi)(x,z)$ of $\Phi$ at $(x,z)$ has the following form:
\[
J(\Phi)(x,z)=\begin{bmatrix}
J\iFop(x)&O_{p,i}\\
\ast&\begin{matrix}
\frac{\partial \widetilde{M}_{n-i+1}}{\partial Z_{n-i+1}}(x,z)&\cdots&\frac{\partial \widetilde{M}_{n-i+1}}{\partial Z_n}(x,z)\\
\vdots&\cdots&\vdots\\
\frac{\partial \widetilde{M}_n}{\partial Z_{n-i+1}}(x,z)&\cdots&\frac{\partial \widetilde{M}_n}{\partial Z_n}(x,z)
\end{matrix}
\end{bmatrix}=
\]
\[
=\begin{bmatrix}
J\iFop(x)&O_{p,i}\\
\ast&\begin{matrix}
m(x)&0&\cdots&0\\
0&m(x)&\cdots&0\\
\vdots&\vdots&\ddots&\vdots\\
0&0&\cdots&m(x)
\end{matrix}\end{bmatrix}
\]
and is of maximal rank $p+i$, since $x$ belongs to $\A^n_m$. In particular $0\in \A^{p+i}$
is a regular value of $\Phi$. From the Weak Transversality Theorem of Thom--Sard we deduce now that there exists a non--empty Zariski open subset $\widetilde{O}$ of $\A^i$ such that for any point $z\in \widetilde{O}$ the equations
\begin{equation}\label{eq:5}
F_1(X)=0 \klk F_p(X)=0, \;\widetilde{M}_{n-i+1}(X,z)=0\klk \widetilde{M}_n(X,z)=0
\end{equation}
intersect transversally at any of their common solutions in $\A^n_m$.\spar

Let $O$ be the image of $\widetilde{O}$ under the given isomorphism which maps $\A^i$ onto $E$. Then $O$ is a non--empty Zariski open subset of $E$.\spar

Fix for the moment a point $z\in \widetilde{O}$. From our previous argumentation and \cite{bank2}, Lemma 1 we conclude that
at any point $x\in \A^n_m$ with
\[
F_1(x)=0 \klk F_p(x)=0, \;\widetilde{M}_{n-i+1}(x,z)=0\klk \widetilde{M}_n(x,z)=0
\]
all $(n-i+1)$--minors of $N(x,z)$ are vanishing.\spar

This means that the equations (\ref{eq:5}) define the intersection of the dual polar variety 
$W_{\ol{K}(\ol{B(z)})}(S)$ with $\A^n_m$. 
Therefore $W_{\ol{K}(\ol{B(z)})}(S)_m$
is either empty or a smooth affine subvariety of $S$ of pure codimension $i$.\spar

Re--parametrizing the algebraic family of affine varieties $\{W_{\ol{K}(\ol{B(z)})}(S)_m\}_{z\in \widetilde{O}}$ by the non--empty Zariski open subset $O$ of $E$ we obtain finally a meagerly generic algebraic family of empty or smooth Zariski open subsets of polar varieties of $S$.\spar

Let $\gamma=(\gamma_1\klk \gamma_n)$ be a generically chosen point of $\Q^n$
and assume that the real variety $S_{\R}$ is smooth. Then we deduce from \cite{bank3} and
\cite{bank4}, Proposition 2 that the generic polar variety $W_{\ol{K}(\ol{\gamma})}(S_{\R})$ is not empty. Hence, without loss of generality, we may assume
$W_{\ol{K}(\ol{\gamma})}(S)_{m}\neq \emptyset$.\spar

For $1\le i \le n-p$ let $z_i:=(\gamma_{n-i+1}\klk \gamma_n)$. Then, similarly as in Example 1, we may argue that
\begin{equation}\label{eq:6}
W_{\ol{K}(\ol{B(z_1)})}(S)_{m}\supset \cdots \supset W_{\ol{K}(\ol{B(z_{n-p})})}(S)_{m}=W_{\ol{K}(\gamma)}(S)_{m}
\end{equation} 
is a descending chain of smooth affine varieties which are codimension one subvarieties one of the other. \spar

Again, in order to retrieve finitely many real algebraic sample points for the connected components of the real variety $S_{\R}$, we may proceed like in Example 1, applying
the general algorithm of \cite{bank2, bank3, bank4} as follows: For any choice of indices $1\le h_1 <\cdots < h_p\le n$ we generate the equations and the inequations that define the descending chain of affine varieties corresponding to (\ref{eq:6}) in order to find the (finitely many) real points contained in the last one. The set of points obtained in this way is again a set of sample points for the connected components of $S_{\R}$.\spar 

A variant of Example 2 is used in \cite{bank5} in order to find
efficiently smooth algebraic sample points for the (non--degenerated)
connected components of singular real hypersurfaces.
\spar

\section{A computational test}\label{s:5}

In this short section we try to get some feeling about the sharpness of the bound 
contained in Proposition \ref{p:a}. To this end we performed a series of computer experiments with generic classic polar varieties of smooth complete intersection manifolds.\spar

More precisely, for a given a triple $(n,p,i)$ of integer parameters with
$1 \le p \le n-1$ and $1 \le i \le n-p$, we did the following:\\
we
\begin{itemize}
\item chose random polynomials $F_1\klk F_p$ of degree two in $\Z[X]$,
\item verified that $F_1\klk F_p$ form a regular sequence in $\Q[X]$
whose set of common zeros $S$ contains no singular point,
\item chose a random integer $((n-p)\times n)$--matrix $a:=\left[
a_{k,l}\right]_{{1 \le k \le n-p-i+1} \atop {1 \le l \le n}}$, and \item
determined the dimension of the singular locus of the classic polar
variety $W_{\ul{K}(\ul{a}^{(i)})}(S)$.
\end{itemize}
The results of such experiments should be interpreted with caution, since
we cannot guarantee that the matrix $a$ satisfies the necessary
genericity condition. However, running the procedure for several random
choices of $a$ may increase our confidence into the experimental results.
\spar

Moreover, the aim was not to test the efficiency of the root finding procedures 
proposed in \cite{bank2, bank3, bank4}, but to check experimentally
a mathematical thesis (the sharpness of the bound contained in Proposition \ref{p:a}). 
Therefore we chose a rather modest sample of equation systems and relied on the most
comfortable software, disregarding its efficiency. \spar

More specifically, 
we chose equations of degree two in order to avoid a complexity explosion
of our computations. In the same vein, we had to control the bit--size of
the coefficients of the polynomials created during our procedure. To get
rid of this situation, we performed our computations modulo the prime
number $q:=10000000019$, taking into account that for any system
$F_1\klk F_p$ of $\Z[X]$ and any integer $((n-p)\times n)$--matrix $a$
there are only finitely many primes $r$ for which the dimension of the
singular locus of $W_{\ul{K}(\ul{a}^{(i)})}(S)$
differs from that of its counterpart obtained over the finite field $\Z_r$
by reducing
$F_1\klk F_p$ and $a$ modulo $r$. Hence we guess that the experimental
results given here reflect the expected generic behavior of the system
$F_1\klk F_p$ over $\Q$.\spar

The computations were performed using the MAGMA package,
on a Core2 DUO processor with 4 Gb of RAM. The search space consisted of
all values of the triple $(n,p,i)$ with $2 \le n \le 11$, $1 \le p \le
n-1$ and $1 \le i \le n-p$. Nevertheless, we had to discard all triples of
the form $(10,p,i)$ with $p=6\klk 9$ and $(11,p,i)$ with $p=5\klk 10$,
since computations ran out of memory.

The experimentation showed that in the hypersurface case, i.e., $p=1$, the
resulting polar varieties were always smooth, as expected.  For $p>1$ we
found that the dimension of the singular locus of the polar varieties was
always $\max\{-1,n-p-(2i+2)\}$, where the dimension of the empty set is
defined by $-1$. In other words, in case $(2i+2)> n-p$, the singular locus
was empty, as predicted for generic polar varieties in Proposition \ref{p:a}, and in case
$(2i+2)< n-p$, the experimentation returned exactly the value $n-p-(2i+2)$
as the dimension of the singular locus, as predicted by the vector bundle viewpoint for
generic classic polar varieties.\spar 

It turned out that in case $2i+2\le n-p$ the singular locus of the polar variety under consideration contained always a real subvariety of the  codimension $2i+2$.\spar

Observe that by Lemma \ref{l:e} 
the value $2i+2$ coincides with the expected codimension of the
stratum $\Sigma^{(i+1)}(\pi_{a^{(i)}})$, in case $2i+2\le n-p$. This finding
suggests that in this case $\bigcup_{i< j\le n-p}\Sigma^{(j)}(\pi_{a^{(i)}}))$
contains the singular locus of $W_{\ul{K}(\ul{a}^{(i)})}\cap S_{reg}$.\vspace{0.4cm}

\textbf{Acknowledgements}\\
The authors are 	much obliged  to two unknown referees who helped by their  comments and hints notably to improve this article.

\end{document}